\documentclass[10pt]{article}
\usepackage{amsmath,amssymb,amsthm,graphicx,tabls}
\usepackage{hyperref}

\newtheorem{thm}{Theorem}[section]

\newtheorem{lem}[thm]{Lemma}
\newtheorem{prop}[thm]{Proposition}

\theoremstyle{definition}
\newtheorem{rem}[thm]{Remark}

\newcommand{\disp}{\displaystyle}

\newcommand{\pa}[1]{\left(#1\right)}
\newcommand{\ac}[1]{\left\{#1\right\}}

\newcommand{\esp}[1]{\mathbb{E}\left[#1\right]}
\newcommand{\pr}[1]{\mathbb{P}\left(#1\right)}

\newcommand{\gep}{\varepsilon}

\newcommand{\pfrac}[1]{\mbox{\small $\frac#1$}}
\newcommand{\To}{\longrightarrow}
\newcommand{\dd}{\textnormal{d}}
\newcommand{\one}{1\hspace{-0.9mm}{\rm l}}
\newcommand{\oc}{[\hspace{-1.3mm}{[}}
\newcommand{\fc}{]\hspace{-1.3mm}{]}}
\newcommand{\Xn}{\widetilde{X}^{(n)}}
\newcommand{\Bn}{B^{(n)}}
\newcommand{\An}{A^{(n)}}
\newcommand{\Mn}{M^{(n)}}

\newcommand{\X}{\widetilde{X}}

\begin{document}

%%%%%%%%%%%%%%%%%%%%%%%%%%%%%%%%%%%%%%%%%%%%%%%%%%%%%%%%%%%%%%%%%%%%%%%%%%
%%%%%%%%%%%%%%%%%%%%%%%%%%%%%%%%%%%%%%%%%%%%%%%%%%%%%%%%%%%%%%%%%%%%%%%%%%

\title{Asymptotic behaviour  of watermelons}
\author{Florent Gillet\thanks{%
Institut \'{E}lie Cartan, Universit\'e Henri Poincar\'e-Nancy I, BP 239,
54506 Vand{\oe }uvre-l\`es-Nancy, France. E-mail: \texttt{%
Florent.Gillet@iecn.u-nancy.fr}}}
\date{}

\maketitle

\begin{abstract}

A watermelon is a set of $p$  Bernoulli paths starting and ending at the same
ordinate, that do not intersect.
In this paper, we show the convergence in distribution of two
sorts of watermelons (with or without wall condition) to processes which
generalize the Brownian bridge and the Brownian excursion in
$\mathbb{R}^p$. These limit processes are defined
by stochastic differential equations. The distributions involved
are those of eigenvalues of some Hermitian random matrices. We give also some
properties of these limit processes.

\end{abstract}

%%%%%%%%%%%%%%%%%%%%%%%%%%%%%%%%%%%%%%%%%%%%%%%%%%%%%%%%%%%%%%%%%%%%%%%%%%
%%%%%%%%%%%%%%%%%%%%%%%%%%%%%%%%%%%%%%%%%%%%%%%%%%%%%%%%%%%%%%%%%%%%%%%%%%

\section{Introduction}
We call $(p,n)$-{\it watermelon} a set of $p$ Bernoulli paths of length $n$
that
meet two conditions. The first condition concerns the starting
points and the endpoints: the $i$-th path starts at level $2i-2$ and
ends at level $k+2i-2$. The integer $k$ is called the {\it deviation} of the
watermelon. The second condition is that the $p$ paths do not touch each
other. An additional condition, called {\it wall condition}, can be
imposed: the
paths should not cross the $x$ axis. Watermelons constitute a particular
configuration of {\it vicious walkers} describing  the situation in which
two or more walkers arriving in the same lattice site annihilate one
another. This model was introduced by Fisher \cite{fisher}, who also
gave a number of physical applications for watermelons.

We shall consider two sorts of $(p,2n)$-watermelons {\it without}
 deviation: {\it
watermelons with wall condition} and {\it watermelons without wall
condition}. They generalize the notion of Dyck path (or excursion) and Grand
Dyck path (or bridge). These processes have many applications, for instance to
    random matrices and  Young tableaux \cite{baik,grabiner,metha}, to the
study of graphs via bijections \cite{bonichon}
and of course to physics \cite{essam,fisher}.

\begin{figure}
\begin{center}
\begin{tabular}{cc}
\includegraphics[width=6cm]{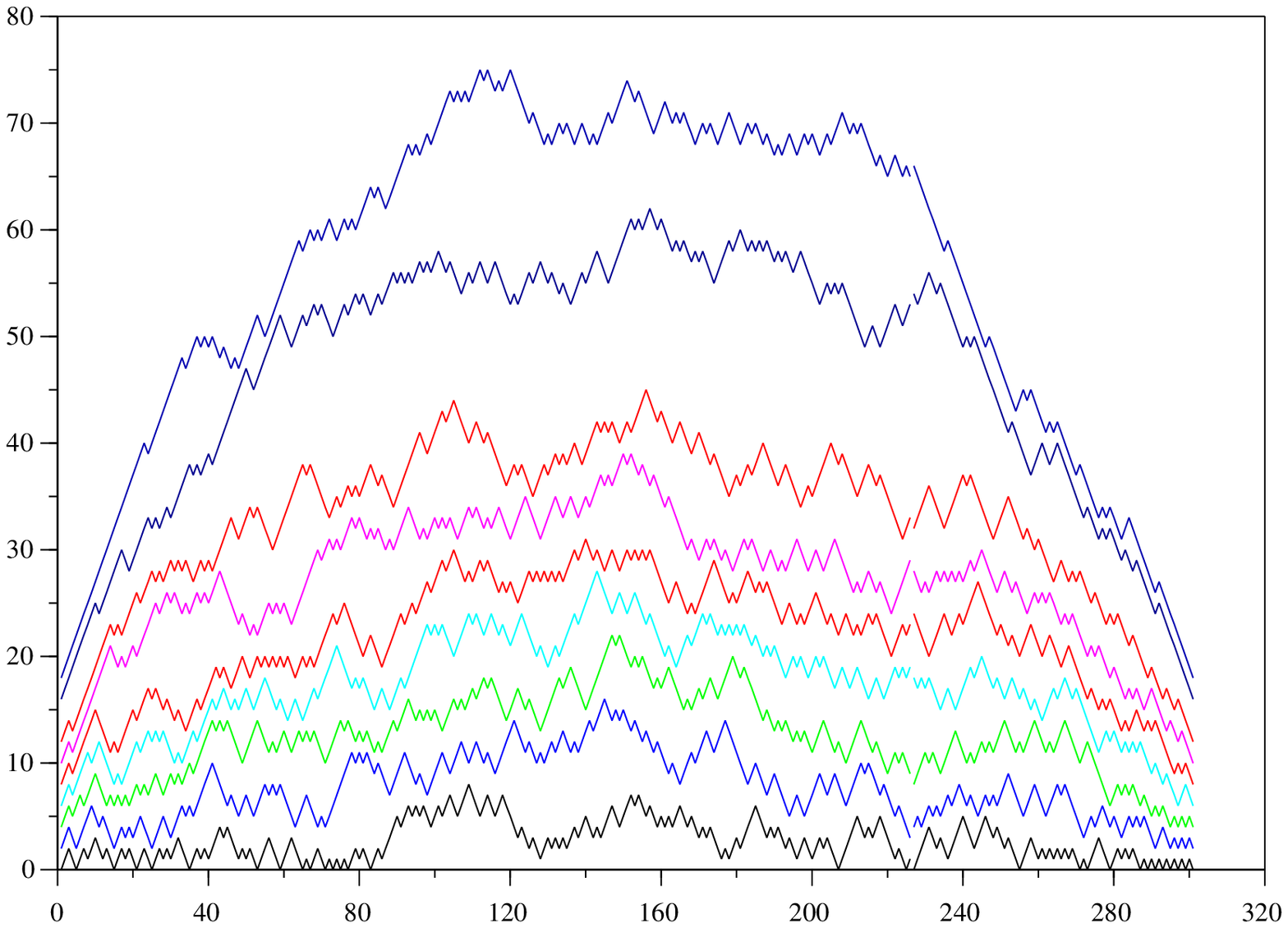}&
\includegraphics[width=6cm]{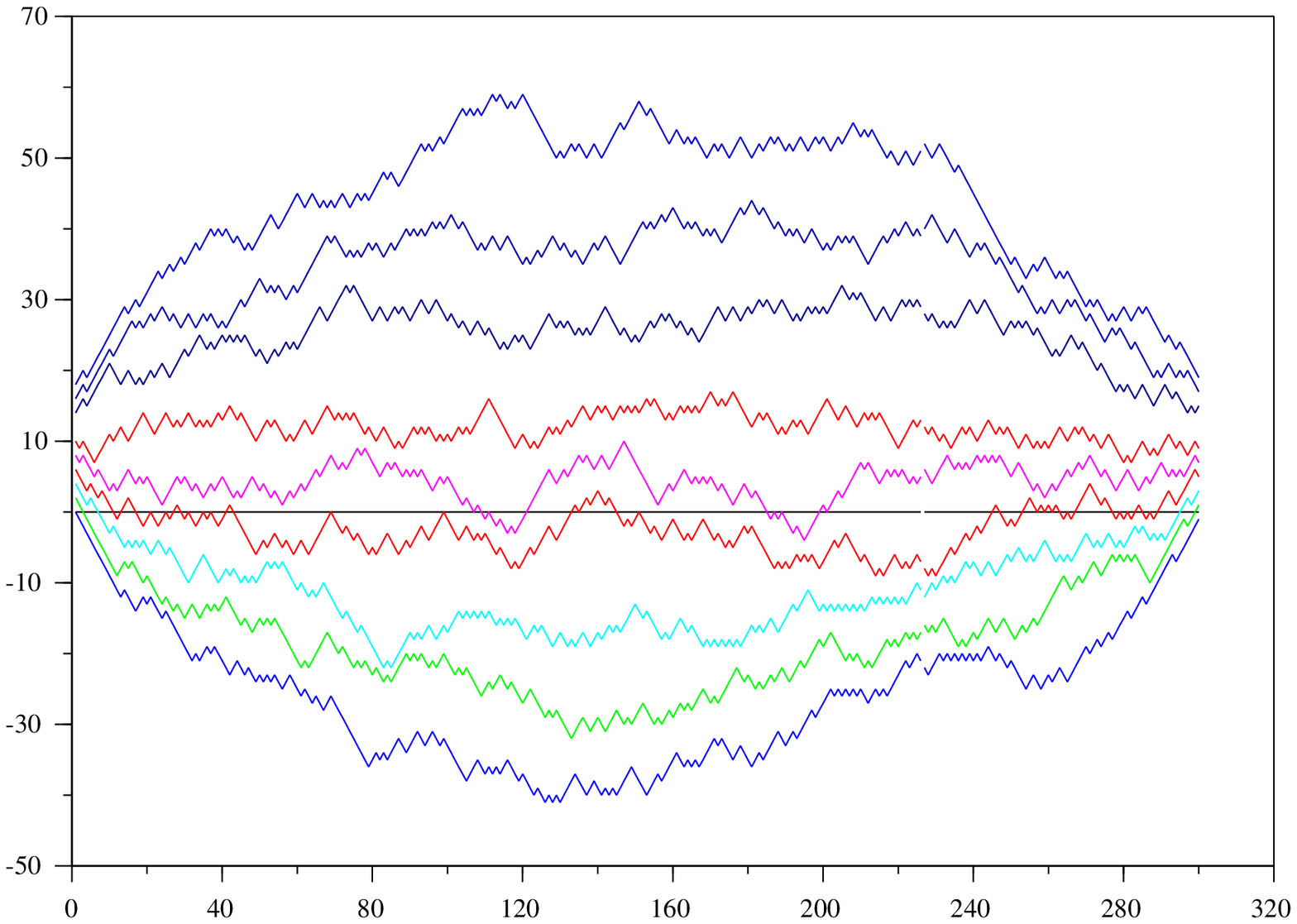}\\
(a)&(b)
\end{tabular}
\caption{$(10,300)$-Watermelons with (a), and without
    (b) wall condition,}
without deviation.
\end{center}
\end{figure}

This paper contains   three sections:
the first section deals with  watermelons with wall condition. We
show that a $(p,2n)$-watermelon   with wall
condition, suitably rescaled in space (by $1/\sqrt{2n}$) and in time
(by $1/(2n)$), converges in distribution to a
stochastic process, that can be seen as
    the analog of a
Brownian excursion for Dyson's version of the  $p$-dimensional  Brownian motion
\cite{dyson}.
This limit process is defined by a stochastic differential
equation. It inherits the  properties of watermelons with wall condition,
namely its branches are  positive (but at $0$ and $1$), and  do not
touch each other. For this  reason  we call the limit process {\it
     continuous $p$-watermelon} with wall condition. We also show, in agreement
with \cite{grabiner},
%\PC{should check that \cite{grabiner} is the good reference}
    that the norm of such continuous $p$-watermelons is a
Bessel bridge with dimension $p(2p+1)$.

In Section \ref{without}, we state similar results for
watermelons without wall condition. The limit
stochastic process is the $p$-dimensional analog of a
Brownian bridge: the intrinsic Brownian bridge of the Weyl chamber
$\{x\in\mathbb{R},x_1<\cdots<x_p\}$ (c.f. \cite{bougerol}).
The norm of a {\it
continuous $p$-watermelon} without wall condition is also a Bessel bridge,
but, as expected,
with dimension $p^2$. The proofs are given in Sections
    \ref{preuve} and  \ref{preuvesm}. They  make
essential use of a variant of Karlin Mac Gregor' formulas
(\cite{karlin,karlinmcgregor,viennot})
for the enumeration of non colliding paths:
 stars with fixed end points.

In the last Section, we give
     the expectation of the elementary symmetric polynomials
of the positions of the continuous $p$-watermelon at time $t$,
    with or without wall
condition. When $p=2$,
we give the exact value of the moments of the two branches at
time $t$.

In \cite{katori}, Katori and Tanemura have shown that the vicious walkers
with and without wall condition converge to the non-colliding
Brownian 
motion (Dyson's Brownian motion \cite{dyson}) and the non-colliding Brownian
meander. In these cases, the one-dimensional distribution of the limit
processes are the eigenvalues of random matrices from a Gaussian Unitary
Ensemble (GUE) and 
from a Gaussian Orthogonal Ensemble (GOE) (\cite{metha}). The watermelons
that we study in this paper can be seen as vicious walkers constrained to
return to $0$ at the end. Once rescaled, our one-dimensional distribusions
are those of the eigenvalues of a Gaussian antisymmetric Hermitian matrices
and a matrices from a GUE.

\section{Main results}
\subsection{Watermelons with wall condition}
\label{with}

Throughout the paper, $W^{(n)}_i(t)$ denotes the 
position at time $t\in[0,2n]$ of the
$i$-th branch of a $(p,2n)$-watermelon  with wall condition, and
$$X^{(n)}_i(t)=\frac1{\sqrt{2n}}W^{(n)}_i([2nt]).$$
We set
$$W^{(n)}(t)=(W^{(n)}_i(t))_{1\leq i\leq p}.$$
Thus
$$W^{(n)}=(W^{(n)}(t)),{t\in[0,2n]})$$
denotes a $(p,2n)$-watermelon  with wall condition that starts
and ends at $0$. We see it as a process in $\mathbb{R}^p$, more
precisely in the Weyl chamber $\{x\in\mathbb{R}^p, 0<
x_1<\cdots<x_p\}$. Similarly
$X^{(n)}$
and $X^{(n)}(t)$
are related to a \textit{renormalized}  $(p,n)$-watermelon with wall condition.

%We denote by
%$\oc x,y \fc$ the set of integers in $[x,y]$.

We give a theorem about a stochastic differential
equation (SDE) which allows us to define properly the limit 
process of $(p,2n)$-watermelons and we state a convergence theorem for
watermelons. A few properties of this limit will be given in the sequel.
%In the second section, we prove the convergence theorem.
These results are proved in section \ref{preuve}.

%\subsection{The results}

Let $H$ be a $(2p+1)\times(2p+1)$ Gaussian antisymmetric Hermitian matrix,
and 
$$\Lambda=(\lambda_1,\dots,\lambda_p)$$
 the $p$ positive eigenvalues of
$H$ sorted by increasing order. We know (cf. \cite{metha}) that the random
variable $\Lambda$ has the density $f$ defined by
%\begin{equation*}
%f(t;x)=\frac{c_p}{\big(t(1-t)\big)^{p^2+p/2}} \prod_{1\leq i<j\leq p}({x_j}%
%^2-{x_i}^2)^2\,\prod_{i=1}^px_i^2 \ \mathrm{e}^{-\frac{\|x\|^2}{2t(1-t)}}1%
%\hspace{-0.9mm}\mathrm{l}_{0\leq x_1\leq \cdots\leq x_p},
%\end{equation*}
\begin{equation*}
f(x)=c_p\prod_{1\leq i<j\leq p}({x_j}%
^2-{x_i}^2)^2\,\prod_{i=1}^px_i^2 \ \mathrm{e}^{-\|x\|^2}1%
\hspace{-0.9mm}\mathrm{l}_{0\leq x_1\leq \cdots\leq x_p},
\end{equation*}
with
\begin{equation*}
c_p=\frac{2^{3p/2}p!}{(2p)!\pi^{p/2} \displaystyle\prod_{1\leq i<j\leq
p}(j-i)(j+i-1)}.
\end{equation*}
We have

\begin{thm}\label{th a p branches}
The stochastic differential equation (SDE)
\begin{align}
\left\{
\begin{array}{l}
\displaystyle \mathnormal{d} X_i(t)=\Big(\frac{-X_i(t)}{1-t}+\frac{1}{X_i(t)}%
+2X_i(t) \sum_{\substack{ j=1  \\ j\neq i}}^{p}\frac{1}{X_i(t)^2-{X_j(t)^2}}%
\Big)\mathnormal{d} s +\mathnormal{d} B_i(t)  \\
X_i(0)=0,\qquad \forall t\in\left]0;1\right[,\ %
X(t)\sim 2\sqrt{t(1-t)}\ \Lambda
\end{array}\hspace{-.4cm}\tag*{$(E_w)$}\label{EDS}
\right.
\end{align}
has a unique strong solution $\big((X_i(t))_{i=1\dots p},t\in\left[0,1%
\right]\big)$.
\end{thm}

\begin{thm}
\label{convergence}
The process $X^{(n)}$ converges in
distribution to the unique 
solution $X$ of the SDE \ref{EDS}.
\end{thm}

This theorem will be proved in the section \ref{demconv}.

\begin{rem}
The three terms of the drift coefficient in \ref{EDS} reflect
the different properties of the watermelons:

\smallskip

-- the first term $-X_i(t)/(1-t)$ pushes the i-th path of $X$
    to $0$ when $t$ is close to $1$.

-- the second term $1/X_i(t)$
keeps the i-th path away from the x-axis.

-- the third term $2X_i(t)/(X_i^2(t)-X_j^2(t))$
keeps the paths away from each other.

The first two terms are also present in the SDE for the normalized Brownian
excursion,
\[
\mathnormal{d} e(t)=\Big(\frac{-e(t)}{1-t}+\frac{1}{e(t)}\Big)\mathnormal{d} s
  +\mathnormal{d} B(t),
\]
while the third is specific to watermelons.
\end{rem}

\begin{rem}
We shall see in section \ref{prop} (Proposition \ref{bessel}) that the
Euclidean norm of the solution of \ref{EDS} is a Bessel bridge (\cite{R-Y})
with dimension $p(2p+1)$.
\end{rem}

\begin{rem}
A solution of \ref{EDS} has the same properties as a watermelon with wall
condition: their paths do not touch each other, stay
positive on $]0,1[$ (cf. Proposition \ref{vander}) also they start and end at
$0$. This is the reason why
we call a solution of \ref{EDS} continuous $p$-watermelon with wall condition.
\end{rem}

%\begin{rem}
%{\bf faire ici le lien entre les watermrlons et les processus de Bougerol}

%\end{rem}

\subsection{Watermelons without wall condition}\label{without}

In this section, we consider $(p,2n)$-watermelon without wall
condition. The results, and the proofs, given in section \ref{preuvesm},
are similar.
Throughout the paper, we denote by
$$\widehat{W}^{(n)}=\big((\widehat{W}_i^{(n)}(t))_{1\leq i\leq
     p},0\leq t\leq 2n\big)$$
a $(p,2n)$-watermelon without wall condition
starting and ending at $0$ and
$$\big(\widehat{X}^{(n)}(t),t\in[0;1]\big)
=\big(\frac1{\sqrt{2n}}\widehat{W}^{(n)}([2nt]))
,t\in[0;1]\big)$$
a renormalized $(p,2n)$-watermelon without wall condition.

Let $\widehat{H}$ be a matrix from a Gaussian Unitary Ensemble, we denote by 
$$\widehat{\Lambda}=(\widehat{\lambda}_1,\dots,\widehat{\lambda}_p)$$
the $p$ eigenvalues of $\widehat{H}$ sorted in increasing order. We know
(cf. \cite{metha}) that the random variable $\widehat{\Lambda}$ has the density 
$\widehat{f}$ defined by
%\begin{equation*}
%\widetilde{f}(t;x)=\frac{2^{-p/2}}{\pi^{p/2}\big(t(1-t)\big)^{p^2/2}
%\prod_{i=1}^{p-1}i!}
%\prod_{1\leq i<j\leq p}(x_j-x_i)^2
%\ \mathrm{e}^{-\frac{\|x\|^2}{2t(1-t)}}1%
%\hspace{-0.9mm}\mathrm{l}_{x_1\leq \cdots\leq x_p}.
%\end{equation*}
\begin{equation*}
\widehat{f}(x)=\frac{2^{-p/2}}{\pi^{p/2}\prod_{i=1}^{p-1}i!}
\prod_{1\leq i<j\leq p}(x_j-x_i)^2
\ \mathrm{e}^{-\|x\|^2}1%
\hspace{-0.9mm}\mathrm{l}_{x_1\leq \cdots\leq x_p}.
\end{equation*}
We have

\begin{thm}\label{th a p branches sm}
The stochastic differential equation
\begin{align}
\left\{
\begin{array}{l}
\displaystyle \mathnormal{d} \widehat{X}_i(t)=\Big(\frac{-\widehat{X}_i(t)}{1-t}%
+\sum_{\substack{ j=1  \\ j\neq i}}^{p}\frac{1}{\widehat{X}_i(t)-{\widehat{X}_j(t)}}%
\Big)\mathnormal{d} s +\mathnormal{d} B_i(t)  \\
\widehat{X}(0)=0,\qquad \forall t\in]0;1[,\ \widehat{X}(t)\sim\sqrt{2t(1-t)}\ \widetilde{\Lambda}
\end{array}\tag*{$(E_{w/o})$}\label{EDSsm}
\right.
\end{align}
has a unique strong solution 
$\widehat{X}=\big((\widehat{X}_i(t))_{i=1\dots p},t\in\left[0,1%
\right]\big)$.

\end{thm}

\begin{thm}\label{convergence sm}

The process $\widehat{X}^{(n)}$ converges in
distribution to the unique  solution  $\widehat{X}$ of  \ref{EDSsm}.
\end{thm}

\begin{rem}
We shall see in section \ref{propsm} (Proposition \ref{bessel sm}) that the
Euclidean norm of the solution of \ref{EDSsm} is a Bessel bridge
with dimension $p^2$.
\end{rem}

\begin{rem}
A solution of \ref{EDSsm} has the same properties as a watermelon without
wall condition:
their paths do not touch each other (cf. Proposition \ref{vander sm}),
start and end at $0$. This is the reason why
we call a solution of \ref{EDSsm} continuous $p$-watermelon without wall condition.
\end{rem}

\begin{rem}
In \cite{bougerol}, Bougerol and Jeulin 
consider the Brownian bridge of
length $T$ on a symmetric space of the non-compact type. They prove that
this rescaled process converges when $T$ tends to infinity, to a limit
process. The generalized radial part of this limit process is the bridge
associated with the intrinsic Brownian motion in a Weyl chamber. For a
suitable choice of the Weyl chamber, a watermelon without wall condition is
the intrinsic Brownian bridge in this chamber. We thus have a geometric 
explanation of the dimension of the Bessel bridge corresponding to the
Euclidean norm of the continuous $p$-watermelon. 
\end{rem}

\section{Watermelons with wall condition: proofs}\label{preuve}
The  existence of solutions of  \ref{EDS} is a consequence of the tightness
of  sequence $X^{(n)}$, and follows from the
proof of  Theorem \ref{convergence}. We focus first on uniqueness
and properties of solutions of  \ref{EDS}.
%\PC{ameliorer ... }

\subsection{Properties}\label{prop}

\begin{prop}
\label{bessel}
Let $X=(X_1,\dots,X_p)$ be a solution to the SDE \ref{EDS}. The
Euclidean norm of
   $X$ is a Bessel bridge
with dimension $p(2p+1)$.
\end{prop}

\begin{proof}%[Proof of Proposition \ref{bessel}]
Set $V=\|X\|$.
Applying It\^o's formula \cite{R-Y} to $V$, we obtain
\begin{align*}
\dd V
&=\sum_{i=1}^p2X_i\dd X_i+\sum_{i=1}^p \dd<X_i,X_i>\\
&=\sum_{i=1}^p\Big(\frac{-2{X_i}^2}{1-t}+2
+\sum_{j=1,j\neq i}^{p}\frac{4{X_i}^2}{{X_i}^2-{X_j}^2}\Big)\dd t
+\sum_{i=1}^p2X_i\dd B_i+p\dd t\\
&=\Big(-\frac{2V}{1-t}+3p
+\sum_{i,j=1,j\neq i}^{p}\frac{4{X_i}^2}{{X_i}^2-{X_j}^2}\Big)\dd t
+\sum_{i=1}^p2X_i\dd B_i.
\end{align*}
Now, we remark that
$$\sum_{i,j=1,j\neq i}^{p}\frac{4{X_i}^2}{{X_i}^2-{X_j}^2}
=2p^2-2p
$$
and that
$$\sum_{i=1}^p2X_i\dd B_i=2\sqrt{V}\dd\beta$$
where $\beta$ is a Brownian motion. We deduce from these equalities that $V$ is
a solution of 
$$
\left\{
\begin{array}{l}
\dd V=\Big(-\frac{2V}{1-t}+p(2p+1)\Big)\dd t+2\sqrt{V}\dd\beta,\\
V_0=0,
\end{array}
\right.$$
while we know that
the squared Bessel bridge with dimension $p(2p+1)$  is another solution
  \cite[Ch. XI]{R-Y}.
\end{proof}

\begin{prop}\label{vander}
If $X=(X_1,\dots,X_p)$ is a solution of  \ref{EDS}, then
$$\pr{\forall t \in \left]0;1\right[,\ 0<X_1(t)<\cdots<X_p(t)<+\infty}=1.$$
\end{prop}
In other words, the paths of a solution of the equation \ref{EDS} do not
touch each other and stay  positive on $\left]0;1\right[$.

\begin{proof}
For every $\gep\in
\left]0;1/2\right[$, we prove  the three  relations:
\begin{align}
&\pr{\forall t \in \left]\gep;1-\gep\right[,\ X_p(t)<+\infty}=1,\label{1}\\
&\pr{\forall t \in \left]\gep;1-\gep\right[,\ X_1(t)<\cdots<X_p(t)}=1,
\label{2}\\
&\pr{\forall t \in \left]\gep;1-\gep\right[,\ X_1(t)>0}=1.\label{3}
\end{align}

Relation (\ref{1}) follows from  Proposition \ref{bessel}:
a Bessel bridge  with dimension $p(2p+1)$ is almost surely finite.

\textit{Proof of relation  }(\ref{2}).
Let $(Z(t))_{\gep\le t\le1-\gep}$ be  defined by
\[
Z(t)=\exp\ac{
\sum_{i=1}^p\int_\gep^t\frac{X_i(s)}{1-s}\dd B_i(s)
-\frac12
\int_\gep^t\sum_{i=1}^p\pa{\frac{X_i(s)}{1-s}}^2\dd s}
\]
and define the probability measure $\mathbb{Q}$ by
$\mathbb{Q}|_{\mathcal{F}_t}=Z(t)\mathbb{P}|_{\mathcal{F}_t}$, i.e.
\[
\forall\Lambda\in\mathcal{F}_t,\ \mathbb{Q}(\Lambda)=\esp{Z(t)\one_\Lambda}.
\]
According to the Girsanov Theorem \cite{R-Y}, under $\mathbb{Q}$,
$X$ is a solution of  an homogeneous SDE:
\[
\dd X_i(t)=
\Big(\frac1{X_i(t)}
+\sum_{\substack{j=1\\j\neq i}}^p\frac{2X_i}{{X_i}^2-{X_j}^2}\Big)dt
+\dd\widetilde{B}(t),
\]
in which $\widetilde{B}$ is a $\mathbb{Q}$-Brownian motion. Set
\[F(x_1,\dots,x_p)=\sum_{1\leq i<j\leq
p}\ln\big({x_j}^2-{x_i}^2\big),\]
  and, for $t\in]\gep,1-\gep[$,
\[U(t)=F(X(t)).\]
Equivalently, relation (\ref{2}) can be written
$$\pr{\forall t \in \left]\gep;1-\gep\right[,\ U(t)>-\infty}=1.$$
It\^o's formula yields that:
\begin{align}\label{ito}
\dd U=
\sum_{k=1}^{p}
\Big(\frac1{X_k(t)}\frac{\partial F}{\partial x_k}
+\big(\frac{\partial F}{\partial x_k}\big)^2
+\frac12\frac{\partial^2 F}{{\partial x_k}^2}\Big)\dd s
+\sum_{k=1}^{p}\frac{\partial F}{\partial x_k}\dd B_k,
\end{align}
in which
\begin{equation*}
\frac{\partial F}{\partial x_k}=\sum_{\substack{j=1\\j\neq k}}^p
\frac{2x_k}{{x_k}^2-{x_j}^2},
\end{equation*}
and
\begin{equation*}
\frac{\partial^2 F}{{\partial x_k}^2}=
\sum_{\substack{j=1\\j\neq k}}^p\Big(\frac{2}{{x_k}^2-{x_j}^2}-
\frac{4{x_k}^2}{({x_k}^2-{x_j}^2)^2}\Big).
\end{equation*}
By symmetry,
\begin{align*}
\sum_{k=1}^{p}\sum_{\substack{j=1\\j\neq k}}^p\frac{1}{{x_k}^2-{x_j}^2}=0
\end{align*}
and
\begin{align*}
\sum_{k=1}^{p}\sum_{\substack{i,j=1\\i\neq j\neq k}}^p
\frac{{x_k}^2}{({x_k}^2-{x_j}^2)({x_k}^2-{x_i}^2)}
=0,
\end{align*}
where $i\neq j\neq k$ means $i\neq j$, $j\neq k$ and $i\neq k$.
This leads to write  (\ref{ito}) under the following form
$$\dd U=S\dd t + M\dd\beta$$
in which  $\beta$ is a $\mathbb{Q}$ Brownian motion,
$$S=2\sum_{\substack{k,j=1\\k\neq j}}^p
\frac{{X_k}^2}{\big({X_k}^2-{X_j}^2\big)^2},$$
and
\[
M
=\sqrt{\sum_{k=1}^{p}\Big(\frac{\partial F}{\partial x_k}\Big)^2}
=\sqrt{2S}.
\]
Let $\tau$ be the stopping time defined by
\begin{align}\label{tau}
\tau &= \inf\left\{\gep\leq t \leq 1-\gep,\ \exists\;i\in\left\oc1,p\right\fc,\
X_i(t)=X_{i+1}(t)\right\}\\
&=\inf\left\{\gep\leq t \leq 1-\gep,\ U(t)=-\infty\right\}\nonumber
\end{align}
with the convention $\tau=1-\gep$ if the infimum is not reach.
We assume that $\mathbb{Q}(\tau<1-\gep)>0$.
As $S$ is positive, by standard  comparison theorems
(see for instance \cite[p. 293]{karatzas}), there exists a solution
$\widetilde{U}$ on $[\gep;1-\gep]$ of the SDE 
\begin{align}\label{ito2}
\left\{
\begin{array}{l}
\dd\widetilde{U}=M\dd\beta\\
\widetilde{U}(\gep)=U(\gep).
\end{array}
\right.
\end{align}
such that
$$\forall t\in\left]\gep,\tau\right[,\ \widetilde{U}(t)\leq U(t)
\qquad\mathbb{Q}-a.s.$$
The previous inequality yields
$$\limsup_{t\to\tau} \widetilde{U}(t) \leq \lim_{t\to\tau} U(t) =-\infty\qquad
\mathbb{Q}\text{-a.s.} $$
so that $\lim_{t\to\tau}\widetilde{U}(t)=-\infty$.
By Dambis Dubins--Schwarz theorem \cite[p.182]{R-Y},
up to an enlargement of the filtered probability space,
there exists a Brownian motion $\widetilde{\beta}$
such that for $t\in\left[\gep;\tau\right[$,
$$\widetilde{U}(t)=\widetilde{U}(\gep)+
\widetilde{\beta}_{\int_\gep^tM(s)\dd s}.$$
Since $\widetilde{U}(t)$ has a limit when $t$ tends to $\tau$, it should be
the same for the time-changed Brownian motion {\it i.e.}
\begin{align*}
0<\mathbb{Q}\big(\tau<1-\gep\big)
\leq\mathbb{Q}\big(\lim_{t\to\tau\wedge1-\gep}U(t)=-\infty\big)
=\mathbb{Q}\big(
\lim_{t\to\tau\wedge1-\gep}\widetilde{\beta}_{\int_\gep^tM(s)}=\infty\big).
\end{align*}
This is absurd
thus $\tau>1-\gep$ $\mathbb{Q}$-a.s. Now $\{\tau > 1-\gep\}$ is
$\mathcal{F}_{1-\gep}$-measurable, hence
$$\mathbb{Q}(\{\tau > 1-\gep\})=\esp{Z_{1-\gep}\one_{\tau >1-\gep}}=1.$$
Since $\esp{Z_{1-\gep}}=1$, we have
$$%\pr{\{\tau \leq 1-\gep\}}=0
%\qquad\text{\it i.e.}\qquad
\pr{\{\tau \geq 1-\gep\}}=1.$$
This evaluation proves the equality (\ref{2}).

\medskip

{\it Proof of relation (\ref{3}).}
Set
\[V=\ln(\prod_{i=1}^pX_i)\]
and
\[
A=\sum_{i=1}^p\frac{1}{{X_i}^2}.
\]
By It\^o's formula, we obtain
\[
\dd V=\left[\frac{-p}{1-t}+\frac A2\right]\dd t
+\sqrt{A}\,\dd \widetilde{B},
\]
where $\widetilde{B}$ is a standard linear Brownian motion.
Set
\begin{align}\label{lowerbound}
\widetilde{V}(t)=V(t)-V(\gep)+\int_\gep^t
\left(\frac{p}{1-s}-\frac A2\right)\dd s,
\end{align}
and
\begin{align}\label{sigma}
\sigma=\inf\big\{t\geq\gep,\ X_1(t)=0\big\}.
\end{align}
Relation (\ref{lowerbound}) has two consequences:
  $\widetilde{V}$ is a time-changed Brownian motion,
\[
\widetilde{V}(t)=
\widetilde{B}\Big(\int_\gep^t\sqrt{A}\ \dd s\Big),
\]
and, if $\sigma\leq1-\gep$,
\[
\lim_{t\to\sigma}\widetilde{V}(t)=-\infty.
\]
Finally,  this last event has null probability,
by standard properties  of Brownian paths,
  entailing
$\pr{\sigma\leq1-\gep}=0$. Relation (\ref{3}) follows.
\end{proof}

\begin{proof}[Proof of  Theorem \ref{th a p branches}]
We shall see in Section \ref{demconv} that the existence of a solution of the SDE
\ref{EDS} is  a consequence of the proof of
Theorem \ref{ethier}.
For $\gep\in\left]0;1/2\right[$, we consider the SDE $(E_\gep)$ defined for
$t\in\left[\gep;1-\gep\right]$ by
\begin{align}
\left\{
\begin{array}{l}
\disp \dd X_i(t)=\Big(\frac{-X_i(t)}{1-t}+\frac{1}{X_i(t)}+2X_i(t)
\sum_{\substack{j=1\\j\neq i}}^{p}\frac{1}{X_i(t)^2-{X_j(t)^2}}\Big) \dd s
+\dd B^i_t \\
X_i(\gep)\sim f(\gep,x) .
\end{array}
\right.
\tag*{($E_\gep$)}\label{EDSe}
\end{align}
It is clear that every solution of \ref{EDS} is a solution of \ref{EDSe}.
The coefficients of the SDE \ref{EDSe} are locally Lipschitz in
$\{0<x_1<\dots<x_p\}$ thus the strong uniqueness holds for the equation
\ref{EDSe} \cite[p. 287]{karatzas} up to the explosion time of the
solution {\it i.e.} on $\left[\gep,\tau\wedge\sigma\right]$ where $\tau$
and $\sigma$ are defined by (\ref{tau}) and (\ref{sigma}). By
Proposition \ref{vander},
    $\tau\wedge\sigma$ is larger than $1-\gep$, thus the
strong uniqueness holds for the equation \ref{EDSe} on
$\left[\gep,1-\gep\right]$.

\smallskip

Let $X$ and $Y$ be two solutions of \ref{EDS}. For every $\gep$ positive,
$X$ and $Y$ are solutions of \ref{EDSe}. Using the uniqueness of
\ref{EDSe}, we have
$$\forall\gep>0,\qquad (X(t),\gep\leq t\leq1-\gep)
\overset{d}{=}  (Y(t),\gep\leq t\leq1-\gep).$$
Hence
$$(X(t),0<t<1)\overset{d}{=}(Y(t),0<t<1).$$
Finally, as $X(0)=Y(0)=X(1)=Y(1)=0$, the two processes $X$ and $Y$ have the
same distribution.
\end{proof}

\subsection{One-dimensional distribution}

Now, we give a few results useful for the proof of Theorem \ref{convergence}.
First, we prove that for a given $t$, the sequence
$\big(X^{(n)}(t)\big)_{n\ge 1}$ converges in distribution to $X(t)$.
%For two vectors $u$ and $v$ in $\mathbb{R}^p$,  $u<v$ means
%$$\forall i\in\left\oc1;p\right\fc,\qquad u_i<v_i.$$
%Also, set
%$$[u;v]=\{x\in\mathbb{R}^p,
%\forall i\in\left\oc1;p\right\fc,\ u_i<x_i<v_i\}.$$

%We have
\begin{prop}\label{densite}
Let $X^{(n)}$ be a renormalized $(p,2n)$-watermelon with
wall condition and let $t$ be in $[0;1]$. $X^{(n)}(t)$ converges in
distribution to $2\sqrt{t(1-t)}\ \Lambda$. 
\end{prop}

\begin{proof} 
We denote by $f(t;x)$ the density function of $2\sqrt{t(1-t)}\ \Lambda$
{\it i.e.}  
\begin{equation*}
f(t;x)=\frac{c_p}{\big(t(1-t)\big)^{p^2+p/2}} \prod_{1\leq i<j\leq p}({x_j}%
^2-{x_i}^2)^2\,\prod_{i=1}^px_i^2 \ \mathrm{e}^{-\frac{\|x\|^2}{2t(1-t)}}1%
\hspace{-0.9mm}\mathrm{l}_{0\leq x_1\leq \cdots\leq x_p},
\end{equation*}
and
\begin{equation*}
c_p=\frac{2^{3p/2}p!}{(2p)!\pi^{p/2} \displaystyle\prod_{1\leq i<j\leq
p}(j-i)(j+i-1)}.
\end{equation*}
For two vectors $u$ and $v$ in $\mathbb{R}^p$,  $u\prec v$ means
$$\forall i\in\left\oc1;p\right\fc,\qquad u_i<v_i.$$
Also, set
$$[u;v]=\{x\in\mathbb{R}^p,
\forall i\in\left\oc1;p\right\fc,\ u_i<x_i<v_i\}$$
where $\oc x;y \fc$ is the set of integers in $[x;y]$. 
Let $u$ and $v$ be two vectors in $\mathbb{R}^p$ such that
$0 \prec u \prec v$ and let $t$ be in $\left[0;1\right]$.
It suffices to prove that
\begin{equation*}
\mathbb{P}\left(X^{(n)}(t)\in\left[u;v\right]\right)
\underset{n\to\infty}{\longrightarrow}\int_{[u;v]}f(t;x)\mathnormal{d} x.
\end{equation*} 

%Let $u$ and $v$ be two vectors of $\mathbb{R}^p$ such that $u<v$ and
%$t\in\left[0;1\right]$.
We have
$$
\pr{X^{(n)}(t)\in\left[u;v\right]}=
\sum_{\substack{u\leq x\leq
       v\\ x\sqrt{2n}\in\mathbb{N}}}
\pr{W^{(n)}([2nt])=x\sqrt{2n}}.
$$

Following \cite{viennot}, we call {\it star}
a set of $p$ random walks, each of length $2n$, that satisfy the non-crossing
condition and start at $0,2,\dots,2p-2$, respectively. The only
difference between a star and a  watermelon is that, in a star, the
$y$-coordinates of the endpoints are unconstrained.

If we denote by $N(m,e)$ the number of stars of length $m$ with wall condition
which end at $(e_1,\dots,e_2)$, we have (Theorem 6. of \cite{viennot})
%\begin{multline}\label{estim}
%N(m,e)=2^{-p^2+p}
%\prod_{i=1}^p\frac{(e_i+1)(m+2i-2)!}
%{(\pfrac{12}(m+e_i)+p)!(\pfrac{12}(m-e_i)+p-1)!}\\
%\times\prod_{1\leq i<j\leq p}(e_j-e_i)(e_j+e_i+2).
%\end{multline}
\begin{multline}\label{estim}
N(m,e)=2^{-p^2+p}
\prod_{i=1}^p(e_i+1)
\!\!\!\prod_{1\leq i<j\leq p}\!\!(e_j-e_i)(e_j+e_i+2)
\prod_{i=1}^p D(m,\frac{e_i}{\sqrt{2n}},p).
\end{multline}
with
\[
D(m,x_i,p)=\frac{(m+2i-2)!}
{(\pfrac{12}(m+x_i\sqrt{2n})+p)!(\pfrac{12}(m-x_i\sqrt{2n})+p-1)!}\cdot
\]
Set $m=[2nt]$. A $(p,2n)$-watermelon which goes through $e$ at time $t$,
can be cut in two 
stars: one of length $m$ and the other of length $2n-m$, both ending at
$e$. The number of $(p,2n)$-watermelons is the number of stars of
length $2n$ which end at $(0,2,\dots,2p-2)$. Regarding this, we
have for $x$ such that $x\sqrt{2n}\in\mathbb{N}^p$
\begin{align*}
\pr{W^{(n)}([2nt])=x\sqrt{2n}}
&=\frac{N(m,x\sqrt{2n})N(2n-m,x\sqrt{2n})}{N(2n,(2i-2)_i)}\\
&=c_{p,n}\prod_{i=1}^p\frac{D(m,x_i,p)D(2n-m,x_i,p)}
{D(2n,\frac{2i-2}{\sqrt{2n}},+p)}
\end{align*}
where
\[\displaystyle
c_{p,n}=\frac{\displaystyle2^{-p^2+p}
\Big[\prod_{1\leq i<j\leq p}(x_j-x_i)(x_i+x_j+\frac{2}{\sqrt{2n}})2n\Big]^2
\Big[\prod_{i=1}^p(x_i+\frac{1}{\sqrt{2n}})\sqrt{2n}\Big]^2}
{\displaystyle\prod_{1\leq i<j\leq p}4(j-i)(j+i-1)\prod_{i=1}^p(2i-1)}.
\]
%and
%\[
%D(m,x_i,p)=\frac{(m+2i-2)!}
%{(\pfrac{12}(m+x_i\sqrt{2n})+p)!(\pfrac{12}(m-x_i\sqrt{2n})+p-1)!}.
%\]

First we have the following estimate
\begin{align*}
c_{p,n}
%=\frac{\displaystyle 2^{-p^2+p}(4n^2)^{p(p-1)/2}(2n)^p}
%{\displaystyle 4^{p(p-1)/2}\prod_{1\leq i<j\leq p}(j-i)(j+i-1)
%\prod_{i=1}^p(2i-1)}\\
%&\hspace{3.6cm}\times\prod_{1\leq i<j\leq
%p}(e_j-e_i)^2(e_i+e_j)^2\prod_{i=1}^p{e_i}^2
%\Big\{1+\mathcal{O}\Big(\frac{1}{\sqrt{n}}\Big)\Big\}\\
=\frac{2^{-p^2+3p}\ p!\ n^{p^2}}
{(2p)!\disp\prod_{1\leq i<j\leq p}(j-i)(j+i-1)}
\prod_{1\leq i<j\leq p}\!\!\!({e_j}^2-{e_i}^2)^2\prod_{i=1}^p{e_i}^2
\Big\{1+\mathcal{O}\Big(\frac{1}{\sqrt{n}}\Big)\Big\}.
\end{align*}

The lemma below gives an estimate of $D(m,e_i,p)$.

\begin{lem}\label{lemme}
Set $n,a,c,d\in\mathbb{N}$, $b\in\frac1{\sqrt{2n}}\mathbb{N}$ and
$t\in\left]0;1\right]$. For $k=nt\{1+\mathcal{O}(\frac{1}{n})\}$, we have
\[
\frac{(2k+a)!}{(k+b\sqrt{2n}+c)!(k-b\sqrt{2n}+d)!}=
\frac{2^{2k+a}}{\sqrt{\pi}}(nt)^{a-c-d-\pfrac{12}}
\,\textnormal{e}^{-2{b}^2/t}
\Big\{1+\mathcal{O}\Big(\frac{1}{\sqrt{n}}\Big)\Big\}
\]
Furthermore, if $a$, $b$, $c$, $d$ and $t$ are in a compact set, then the
$\mathcal{O}$ is uniform in all these variables.
%\PC{ici $t$ doit etre dans un compact de $(0,1]$ et non pas de $[0,1]$,
%si je ne me trompe.}
\end{lem}

\begin{proof}
Using Stirling Formula and an asymptotic expansion of the logarithm, we obtain
\[
(k+b\sqrt{2n}+c)!
=\sqrt{2\pi k}\ k^{k+b\sqrt{2n}+c}\,\textnormal{e}^{-k+\frac{2b^2}{t}}
\Big\{1+\mathcal{O}\Big(\frac{1}{\sqrt{n}}\Big)\Big\}.
\]
We have in the same way
\[
(k-b\sqrt{2n}+d)!
=\sqrt{2\pi k}\ k^{k-b\sqrt{2n}+d}\,\textnormal{e}^{-k+\frac{2b^2}{t}}
\Big\{1+\mathcal{O}\Big(\frac{1}{\sqrt{n}}\Big)\Big\}
\]
and
\[
(2k+a)!
=\sqrt{4\pi k}\ (2k)^{2k+a}\,\textnormal{e}^{-2k}
\Big\{1+\mathcal{O}\Big(\frac{1}{\sqrt{n}}\Big)\Big\},
\]
leading to the result, after simplification.
\end{proof}

This lemma yields the following estimate of $D(m,x_i,p)$:
\[
D(m,x_i,p)=\frac1{\sqrt{\pi}}(nt)^{2i-2p-\pfrac{32}}2^{m+2i-2}
\,\textnormal{e}^{-{x_i}^2/(2t)}
\Big\{1+\mathcal{O}\Big(\frac{1}{\sqrt{n}}\Big)\Big\}
\]
and so
\begin{align*}
\frac{D(m,x_i,p)D(2n-m,x_i,p)}{D(2n,\frac{2i-2}{\sqrt{2n}},p)}&\\
=\frac1{\sqrt{\pi}}\big(nt(1-t)&\big)^{2i-2p-\pfrac{32}}2^{2i-2}
\ \textnormal{e}^{-{x_i}^2/(2t(1-t))}
\Big\{1+\mathcal{O}\Big(\frac{1}{\sqrt{n}}\Big)\Big\}.
\end{align*}
We deduce  that
\begin{align*}
\mathbb{P}\Big(
W^{(n)}([2nt])=&x\sqrt{2n}\Big)\\
&=c_{p,n}\prod_{i=1}^p\frac{D(m,x_i,p)D(2n-m,x_i,p)}{D(2n,0,i-1+p)}\\
&=c_{p}\frac
{\prod_{1\leq i<j\leq p}({x_j}^2-{x_i}^2)^2\prod_{i=1}^p{x_i}^2}
{(2n)^{p/2}\big(t(1-t)\big)^{p^2+p/2}}
\ \textnormal{e}^{-\frac{\|x\|^2}{2t(1-t)}}
\Big\{1+\mathcal{O}\Big(\frac{1}{\sqrt{n}}\Big)\Big\}\\
&=\Big(\frac{2}{n}\Big)^{p/2}
f(t;x)\Big\{1+\mathcal{O}\Big(\frac{1}{\sqrt{n}}\Big)\Big\}
\end{align*}
and finally
\begin{align*}
\pr{X^{(n)}(t)\in\left[u,v\right]}
=\Big(\frac{2}{n}\Big)^{p/2}\sum_{\substack{u\leq x\leq
       v\\e\sqrt{2n}\in\mathbb{N}}}
f(t;x)\Big\{1+\mathcal{O}\Big(\frac{1}{\sqrt{n}}\Big)\Big\}.
\end{align*}

In the equality above, the $\mathcal{O}$ is uniform in $x$ and can be
factorized. Then, this expression becomes a Riemann sum which tends to
$\int_{[u;v]}f(t;x)\dd x$. Indeed, since the $x_i\sqrt{2n}$ have the same
parity, the number of terms of a volume unit is
$(\frac{\sqrt{2n}}{2})^{p}=(\frac{n}{2})^{p/2}$.
\end{proof}

\subsection{Proof of Theorem \ref{convergence}}\label{demconv}

The following lemma explains that the probability that the branches of a
$(p,2n)$-water\-melon with and wall condition do not near $0$
or are relatively far from each other, tends to $1$ as $n$ tends to infinity.
In the proof of Theorem \ref{convergence}, this lemma allows the restriction to such watermelons.
%\PC{to be improved}

For $\gep\in \left]0,1\right[$ and $\alpha\in \left]0,1/4\right[$, we
define the sets
\begin{align*}
\Lambda_n^1(\alpha)
&=\left\{\exists k\in\oc2\gep n,((1-\gep)\wedge\tau_n^r)2n\fc,\
\exists i\in\oc1,p\fc \textnormal{ s.t. } W_i^{(n)}(k)\leq n^\alpha
\right\}\\
&=\left\{\exists k\in
\oc2\gep n,((1-\gep)\wedge\tau_n^r)2n\fc,
\textnormal{ s.t. } W_1^{(n)}(k)\leq n^\alpha
\right\}
\end{align*}
and
\begin{multline*}
\Lambda_n^2(\alpha)
=\left\{\exists k\in\oc2\gep n,((1-\gep)\wedge\tau_n^r)2n\fc,\
\exists i,j\in\oc1,p\fc\ \textnormal{ s.t. }i<j \text{ and}\right.\\
\hspace{6cm}\left.W_j^{(n)}(k)^2-W_i^{(n)}(k)^2\leq n^\alpha\right\}
\end{multline*}
where $r>0$ and $\tau_n^r=\inf\{t>\gep,\ X^{(n)}_p(t)\leq r\}$ .
We have

\begin{lem}\label{restriction}
For $i\in \{1,2\}$,
\[
\lim_{n\to\infty}\pr{\Lambda_n^i(\alpha)}=0.
\]
\end{lem}

\begin{proof}
We have shown in Proposition \ref{densite} that
\[
\pr{X^{(n)}(t)\in\left[u;v\right]}=
\frac1{(2n)^{p/2}}\sum_{\substack{u\sqrt{2n}\leq x\sqrt{2n}\leq
v\sqrt{2n}\\x\sqrt{2n}\in\mathbb{N}}}f(t;x)
\Big\{1+\mathcal{O}\Big(\frac{1}{\sqrt{n}}\Big)\Big\}.
\]
Thus for $t\in\left[\gep,(1-\gep)\right]$, 
we have
\begin{align*}
\mathbb{P}\Big(W_1^{(n)}([2nt])\in\left[0,n^\alpha\right]&
,\ t\leq\tau_n^r\Big)\\
&\leq\pr{X^{(n)}(t)
\in\left[0,n^{\alpha-1/2}\right]\times\left[0,r\right]^{p-1}}\\
&=\frac{1}{(2n)^{p/2}}
\sum_{\substack{0\leq x_1\leq n^{\alpha-1/2}\\x_1\sqrt{2n}\in\mathbb{N}}}
\sum_{\substack{0\leq x_i\leq r\\x_i\sqrt{2n}\in\mathbb{N}\\i=2\dots p}}
f(t;x)\Big\{1+\mathcal{O}\Big(\frac{1}{\sqrt{n}}\Big)\Big\}.
\end{align*}
Since $t\in[\gep,(1-\gep)]$ and $x_1\leq n^{\alpha-1/2}$,
\begin{align*}
f(t;x)
&\leq c_{p}
\frac
{\prod_{1\leq i<j\leq p}({x_j}^2-{x_i}^2)^2\prod_{i=1}^p{x_i}^2}
{\big(\gep(1-\gep)\big)^{p^2+p/2}}
\ \textnormal{e}^{-2({x_2}^2+\cdots+{x_p}^2)}\\
&\leq \sum_{k=2}^{4p-2}{x_1}^kf_k(x_2,\dots,x_p)
\leq\sum_{k=2}^{4p-2}{x_1}^2f_k(x_2,\dots,x_p)\\
&\leq n^{2\alpha-1}\sum_{k=2}^{4p-2}f_k
\end{align*}
where $f_k$ is the product of a polynomial with positive coefficients and
of an exponential term
$\exp(-2({x_2}^2+\cdots+{x_p}^2))$, in particular $f_k$ is
integrable. Using the previous majoration, we obtain
\begin{multline*}
\pr{W_1^{(n)}([2nt])\in\left[0,n^\alpha\right]
,\ t\leq\tau_n^r}\\
\leq\frac{1}{\sqrt{2n}}
\sum_{\substack{0\leq x_1\leq n^{\alpha-1/2}\\x_1\sqrt{2n}\in\mathbb{N}}}
\frac{1}{(2n)^{(p-1)/2}}
\sum_{\substack{0\leq x_i\leq r\\x_i\sqrt{2n}\in\mathbb{N}\\i=2\dots p}}
n^{2\alpha-1}\sum_{k=2}^{4p-2}f_k
\Big\{1+\mathcal{O}\Big(\frac{1}{\sqrt{n}}\Big)\Big\}.
\end{multline*}
The sum over the $x_i$ for $i\in\{2,\dots,p\}$ is a Riemann sum, thus

\[
\frac{1}{(2n)^{(p-1)/2}}
\sum_{\substack{0\leq x_i\leq r\\x_i\sqrt{2n}\in\mathbb{N}\\i=2\dots p}}
\sum_{k=2}^{4p-2}f_k
=\int_0^r\sum_{k=2}^{4p-2}
f_k(x)\dd x\ \left\{1+\mathcal{O}\left(\frac 1n\right)\right\}
\leq c
\]
where $c$ is a positive constant, leading to:
\[
\pr{W_1^{(n)}([2nt])\in\left[0,n^\alpha\right],\ t\leq\tau_n^r}
\leq c n^{3(\alpha-1/2)}.
\]
Finally, for $\alpha\le 1/4$,
\begin{align*}
\pr{\Lambda_n^1(\alpha)}
&\leq\pr{\exists k\leq2\big((1-\gep)\wedge\tau_n^r-\gep\big)n^{1-\alpha}
,\ W_1^{(n)}(2\gep n+kn^\alpha)\leq2n^\alpha}\\
&\leq \sum_{k=0}^{(1-2\gep)2n^{1-\alpha}}
\pr{W_1^{(n)}(2\gep n +kn^\alpha)\leq 2n^\alpha
,\ \gep+\frac{kn^\alpha}{2n}\leq\tau_n^r}\\
&\leq 2cn^{3(\alpha-1/2)}(1-2\gep)n^{1-\alpha}\\
&\leq 2c n^{2\alpha-1/2}.
\end{align*}
entailing  the first point.

For $i=2$,
we  proceed  as above.
Set
\[
\Lambda_n^2(\alpha)=\bigcup_{j=1}^{p-1}\Lambda_n^{2,j},
\]
in which
\[
\Lambda_n^{2,j}=\left\{\exists k\in\oc2\gep n,((1-\gep)\wedge\tau_n^r)2n\fc,
\text{ s.t. }W_{j+1}^{(n)}(k)^2-W_{j}^{(n)}(k)^2\leq n^\alpha\right\}.
\]
Let us fix $j\in[1,p-1]$.
For $t\in[\gep,(1-\gep)]$, Proposition \ref{densite} yields
\begin{multline*}
\pr{W_{j+1}([2nt])^2-W_{j}([2nt])^2\in[0,n^\alpha],\ t\leq\tau_n^r}\\
=\frac{1}{(2n)^{p/2}}
\sum_{\substack{0\leq x_i\leq
r,\ x_i\sqrt{2n}\in\mathbb{N}\\{x_{j+1}}^2-{x_{j}}^2\leq n^{\alpha-1}/2}}
f(t;x)\Big\{1+\mathcal{O}\Big(\frac{1}{\sqrt{n}}\Big)\Big\}.
\end{multline*}
replacing ${x_{j+1}}^2$ with $({x_{j+1}}^2-{x_{j}}^2)+{x_{j}}^2$ in
$f(t;x)$, we  bound  this function when $t\in[\gep,1-\gep]$
and ${x_{j+1}}^2-{x_{j}}^2\leq n^{\alpha-1/2}$:
\begin{align*}
f(t;x)&\leq c_p\frac
{\prod_{1\leq i<j\leq p}({x_j}^2-{x_i}^2)^2\prod_{i=1}^p{x_i}^2}
{\big(\gep(1-\gep)\big)^{p^2+p/2}}
\ \textnormal{e}^{-2({x_1}^2+\cdots+{x_p}^2)}\\
&\leq \sum_{k=2}^{4p-2}\Big({x_{j+1}}^2-{x_{j}}^2\Big)^k
f_k(x_1,\dots,x_j,x_{j+2},\dots,x_p)
\leq n^{2\alpha-1}\sum_{k=2}^{4p-2}f_k
\end{align*}
where the $f_k$ are  products of an exponential 
term $\exp\big(-2({x_1}^2+\dots+{x_{j}}^2+{x_{j+2}}^2+\dots+{x_{p}}^2)\big)$
with polynomials.
Set $y={x_{j+1}}^2-{x_{j}}^2$, we have
\begin{align*}
&\pr{W_1^{(n)}([2nt])\in\left[0,n^\alpha\right],\ t\leq\tau_n^r}\\
&\leq\frac{1}{\sqrt{2n}}
\sum_{\substack{0\leq x_j\leq r\\x_j\sqrt{2n}\in\mathbb{N}}}
\sum_{\substack{0\leq y\leq n^{\alpha-1}\\2ny\in\mathbb{N}}}
\frac{1}{(2n)^{(p-1)/2}}
\sum_{\substack{0\leq x_i\leq r\\x_i\sqrt{2n}\in\mathbb{N}\\i\neq j,j+1}}
n^{2\alpha-1}\sum_{k=2}^{4p-2}f_k
\Big\{1+\mathcal{O}\Big(\frac{1}{\sqrt{n}}\Big)\Big\}\\
&\leq n^{3(\alpha-1/2)}\frac{1}{(2n)^{(p-1)/2}}
\sum_{\substack{0\leq x_i\leq r\\x_i\sqrt{2n}\in\mathbb{N}\\i\neq j+1}}
\sum_{k=2}^{4p-2}f_k
\Big\{1+\mathcal{O}\Big(\frac{1}{\sqrt{n}}\Big)\Big\}.
\end{align*}
Now the $\mathcal{O}$ are uniform and the sum over the $x_i$ is a Riemann
sum thus
\[
\frac{1}{(2n)^{(p-1)/2}}
\sum_{\substack{0\leq x_i\leq r\\xi\sqrt{2n}\in\mathbb{N}\\i\neq j+1}}
\sum_{k=2}^{4p-2}f_k
=\int_{[0;r]^p}\sum_{k=2}^{4p-2}
f_k(x)\dd x\ \left\{1+\mathcal{O}\left(\frac 1n\right)\right\}
\leq c
\]
where $c$ is a positive constant. This provides the following upper bound
\[
\pr{W_{j+1}^{(n)}([2nt])^2-W_{j}^{(n)}([2nt])^2\in\left[0,n^\alpha\right]
,\ t\leq\tau_n^r}
\leq c n^{3(\alpha-1/2)}.
\]

To complete this proof,  notice that we have
\begin{align*}
\mathbb{P}\Big(&\Lambda_n^{2,i}(\alpha)\Big)\\
&\leq\pr{\exists k\leq2\big((1-\gep)\wedge\tau_n^r-\gep\big)n^{1-\alpha},\ %
W_{j+1}^{(n)}(t_k)^2-W_{j}^{(n)}(t_k)^2
\leq3n^\alpha}\\
&\leq \sum_{k=0}^{(1-2\gep)n^{1-\alpha}}
\pr{W_{j+1}^{(n)}(t_k)^2-W_{j}^{(n)}(t_k)^2
\leq 3n^\alpha, t_k\leq\tau_n^r}\\
&\leq c'n^{3(\alpha-1/2)}(1-2\gep)n^{1-\alpha}\leq c'' n^{2\alpha-1/2}
\end{align*}
where $t_k=2\gep n+kn^\alpha$, $c'$ and $c''$ are constants and that $\alpha$ is less than $1/4$.

\end{proof}
\begin{proof}[Proof of Theorem \ref{convergence}]
Let 
$$\Xn=(\X^{(n)}(t),t\in[0;1])$$
 be the piecewise constant process in
$\mathbb{R}^p$ such that for every
$k\in\oc0,2n\fc$ and every $i\in\oc1,p\fc$,
\[
\Xn_i(\frac{k}{2n})=\frac{1}{\sqrt{2n}}(W_i^{(n)}(k)+1).
\]
We consider here sample paths that belong to the space $\mathcal{D}$ of cad-lag
(right continuous left limit) functions endowed with the Skorokhod topology
(cf. \cite[ch 2]{billingsley}).
To prove this theorem, it suffices to show  that the process
$(\Xn(t))_t$ converges in distribution to the process $X$. The
convergence in distribution of the process will be a consequence of Theorem
\ref{ethier} below.
$(\Xn(t),t\in[\gep,1-\gep])$ to the process $(X(t),t\in[\gep,1-\gep])$
for every $\gep\in]0;1/2[$. First, let us show that this convergence
entails Theorem \ref{convergence}. Then we shall state and prove Theorem
\ref{ethier}. Finally we shall apply this theorem to our case.

If we prove that $\Xn$ converges to
$X$ on $[\gep;1-\gep]$ for every $\gep\in]0;1/2[$, then the uniqueness
of the solution of \ref{EDS} implies that $(\Xn(t),t\in]0;1[)$ converges to
$(X(t),t\in]0;1[)$. Since $X$ tends to $0$ when $t$ tends to $0$ and $1$,
we have the convergence on $[0;1]$. 

\smallskip

Theorem \ref{ethier} is similar to Theorem 4.1 in \cite{ethier}.

\begin{thm}\label{ethier}
We assume the strong uniqueness of the solution of the SDE
\begin{align}  \label{edsethier}
\mathnormal{d} X(t)=b(X(t),t)\mathnormal{d}t+\sigma(X(t),t)\mathnormal{d}B(t).
\end{align}
which satisfies a property $\mathcal{P}$.
Let $X^{(n)}$ and $B^{(n)}$ be two processes with cad-lag sample paths
on $\left[0,1\right]$ and let $A^{(n)}=(A_{i,j}^{(n)})$ be a symmetric
$d\times d$ matrix-valued process such that $A_{i,j}^{(n)}$ has
cad-lag sample paths on $\left[0,1\right]$ and $A^{(n)}(t)-A^{(n)}(s)$ is
nonnegative definite for $t>s\geq 0$. Set
\[\tau_n^r=\inf\{t>0\text{ such that } |X^{(n)}(t)|\geq r \ or\
|X^{(n)}(t^-)|\geq r\}\]
and $\mathcal{F}^n_t=\sigma(X^{(n)}(s),B^{(n)}(s),A^{(n)}(s),s\leq t)$.

We assume that the following assumptions hold for every $r>0$, $T>0$ and
$i,j=1,\dots,d$:

\begin{gather*}
M^{(n)}_i=X_i^{(n)}-B_i^{(n)} \text{ is a $\mathcal{F}_t^n$-local
martingale }  \tag*{(H1)}
\label{h1} \\
M_i^{(n)}M_j^{(n)}-A_{ij}^{(n)}
\text{ is a $\mathcal{F}_t^n$-local martingale}
\tag*{(H2)}  \label{h2} \\
\lim_{n\to\infty} \mathbb{E}\bigg[\sup_{t\leq T\wedge\tau_n^r}\big|%
X^{(n)}(t)-X^{(n)}(t^-)\big|^2\bigg]=0  \tag*{(H3)}  \label{h3} \\
\lim_{n\to\infty} \mathbb{E}\bigg[\sup_{t\leq T\wedge\tau_n^r}\big|%
B^{(n)}(t)-B^{(n)}(t^-)\big|^2\bigg]=0  \tag*{(H4)}  \label{h4} \\
\lim_{n\to\infty} \mathbb{E}\bigg[\sup_{t\leq T\wedge\tau_n^r} \big|%
A_{ij}^{(n)}(t)-A_{ij}^{(n)}(t^-)\big|^2\bigg]=0  \tag*{(H5)}  \label{h5} \\
\sup_{t\leq T\wedge\tau_n^r}\bigg|B_i^{(n)}(t)-\int_0^t
b_i(X^{(n)}(s),s)\mathnormal{d}
s\bigg| \overset{\mathbb{P}}{\longrightarrow}0  \tag*{(H6)}  \label{h6} \\
\sup_{t\leq T\wedge\tau_n^r}\bigg|A_{ij}^{(n)}(t)-\int_0^t
a_{ij}(X^{(n)}(s),s)%
\mathnormal{d} s \bigg|\overset{\mathbb{P}}{\longrightarrow}0  \tag*{(H7)}
\label{h7}
\end{gather*}
where $(a_{ij})=\sigma\sigma^*$. Moreover, we assume that every
accumulation point of $X^{(n)}$ satisfies the property $\mathcal{P}$
and  that $X^{(n)}(0)$ converges in distribution to $X(0)$.

Then the process $(X^{(n)})_n$ converges in distribution to the process $X$.
\end{thm}

\begin{proof}
We shall prove this theorem in two steps whether the coefficients are
homogeneous in time or not.

\medskip

-- Let us assume that the coefficients $\sigma$ and $b$ are homogeneous in
time.

Using the same arguments as in the proof of Theorem 4.1 in \cite{ethier}, the
sequence $\big(X^{(n)}(.\wedge\tau_n^r)\big)$ is relatively compact and a
limit $X^{r_0}$ of a subsequence of $\big(X^{(n)}(.\wedge\tau_n^{r_0})\big)$
is a solution of (\ref{edsethier}).

As $X^{r_0}$ is an accumulation point
of the sequence $X^{(n)}$, the property
$\mathcal{P}$ holds for $X^{r_0}$, thus, by uniqueness, $X^{r_0}$ has the
same distribution as $X$. Therefore $X^{(n)}(.\wedge\tau_n^r)$ converges in
distribution to $X(.\wedge\tau^r)$ for every $r>0$. But $\tau^r$ tends to
infinity along with $r$. Thus $X^{(n)}$ converges in distribution to $X$.

\medskip

-- When the coefficients $\sigma$ and $b$ are not homogeneous.

We consider the time-space process
$\bar{X}^{(n)}(t)=\big(X^{(n)}(t),t\big)$. It is easy to show that if
every hypothesis is satisfied for $X^{(n)}$, $B^{(n)}$, $A^{(n)}$, $b$ and
$\sigma$, then they hold for $\bar{X}_n$, the processes $\bar{B}$
and $\bar{A}^{(n)}$ and the functions $\bar{b}$ et
$\bar{\sigma}$ defined by
\begin{gather*}
\bar{B}^{(n)}=\binom{B^{(n)}}{t},\qquad
\bar{b}\binom{x}{t}=\binom{b(x,t)}{t},\qquad
\bar{\sigma}\binom{x}{t}=\binom{\sigma(x,t)}{0},\\
\bar{A}^{(n)}_{ij}=
\left\{\begin{array}{ll}
A^{(n)}_{ij} &\text{if }i,j\in\oc1,d\fc\\
0 &\text{if }i\text{ or }j=d+1\\
t/n&\text{if }i=j=d+1
\end{array}\right..
\end{gather*}
Hence $\bar{X}^{(n)}$ converges in distribution to $(X,t)$, therefore
$X^{(n)}$ converges in distribution to $X$.
\end{proof}

\medskip

\begin{rem}
Theorem \ref{ethier}
does not only give  the convergence of $X^{(n)}$ to $X$, but also, by
the relatively compactness of  $X^{(n)}$, the existence of a solution of
the SDE \ref{EDS}.
\end{rem}

To prove Theorem \ref{convergence}, it suffices to exhibit some processes $A$,
$B$ and $M$ and verify the assumptions of Theorem \ref{ethier} with the
property $\mathcal{P}=\{$the 1-dimensional distribution of $X$ at time $t$
has the density $f(t;x)$$\}$.
In point of fact, we want to prove the convergence of $\X^{(n)}$ only on
$[\gep,1-\gep]$. Hence, it suffices to verify the assumptions for
$t\geq\gep$ and the assumption \ref{h6''} and \ref{h7''} instead of \ref{h6}
and \ref{h7}.
\begin{gather}
\sup_{\gep\leq t\leq T\wedge\tau_n^r}\bigg|B_i^{(n)}(t)-B_i^{(n)}(\gep)-\int_\gep^t
b_i(\X^{(n)}(s),s)\mathnormal{d}
s\bigg| \overset{\mathbb{P}}{\longrightarrow}0  \tag*{(H6')}  \label{h6''} \\
\sup_{\gep\leq t\leq T\wedge\tau_n^r}\bigg|A_{ij}^{(n)}(t)-A_{ij}^{(n)}(\gep)
-\int_\gep^t a_{ij}(\X^{(n)}(s),s)%
\mathnormal{d} s \bigg|\overset{\mathbb{P}}{\longrightarrow}0  \tag*{(H7')}
\label{h7''}
\end{gather}
To proof this, we apply Theorem \ref{ethier} with $\Xn(\gep+t)$,
$B_i^{(n)}(\gep+t)-B_i^{(n)}(\gep)$ and 
$A_{i,j}^{(n)}(t)-A_{i,j}^{(n)}(\gep)$.
%, and apply
%Theorem \ref{ethier} with $\hat{X}^{(n)}$, $\hat{B}_i^{(n)}$ and
%$\bar{A}_{i,j}^{(n)}$.

First, let us give the definitions of the processes $A$, $B$ and $M$, the
stopping times $T$ and $\tau$ and the functions  $b$, $\sigma$ and $a$: for
$t\in[0;1]$,
\begin{gather*}
\Bn(t)=(\Bn_i(t))_i=\Big(\sum_{l=0}^{[2nt]-1}
\mathbb{E}\left[\Xn_i(\frac{l+1}{2n})-
\Xn_i(\frac{l}{2n})\Big|%\Xn(\frac{l}{2n})\right]\Big)_i,\\
\mathcal{F}_l\right]\Big)_i,\\
\Mn=\Xn-\Bn,\\
\An_{i,j}(t)=\!\!
\sum_{l=0}^{[2nt]-1}\mathbb{E}\left[\Mn_i(\frac{l+1}{2n})\Mn_j(\frac{l+1}{2n})
-\Mn_i(\frac{l}{2n})\Mn_j(\frac{l}{2n})\Big|%\Mn(\frac{l}{2n})\right], \\
\mathcal{F}_l\right], \\
\mathcal{F}_l=\sigma\big(\Xn(\frac{j}{2n}),j=0,\dots,l\big)\\
\text{for $r>0$ and $\varepsilon>0$,}\qquad\tau_n^r=\inf_{t\geq\gep}%
\{\Xn_p(t)\geq r\}, \qquad T=1-\varepsilon, \\
b(t,x)=(b_i(t,x))_{1\leq i\leq p} =\Big(\frac{-x_i}{1-t}+\frac{1}{x_i}%
+\sum_{j=1}^p\frac{2x_i}{x_i^2-x_j^2}\Big) _{1\leq i\leq p}\\
\text{and}%
\qquad \sigma(x,y,t)=a(x,y,t)=Id_2.
\end{gather*}

Some assumptions of Theorem \ref{ethier} are easily verified. First, it is
clear, by definition, that $\Mn$ and $\Mn_i\Mn_j-\An_{i,j}$ are
both martingales (\ref{h1} and \ref{h2}). After, the assumptions \ref{h3},
\ref{h4} and \ref{h5} are easy to state. Indeed, the paths of processes
$\Xn_i$ are piecewise continuous functions with jumps of amplitude
$1/\sqrt{2n}$. It is the same for the paths of $\Bn_i$, and, for the
processes $\An_{i,j}$, it is easy to show that on $[0,\tau_n^r]$, we have
$|\An_{i,j}((k+1)/2n)-\An_{i,j}(k/2n)|\leq3r/\sqrt{2n}$. Hence the three
supremums of these assumptions are bounded above by a
$\mathcal{O}\big(1/\sqrt{2n}\big)$. Finally, Proposition \ref{densite}
prove that $\Xn(\gep)$ tends in distribution to $X(\gep)$ and that
every accumulation point of $\Xn$ verifies the property $\mathcal{P}$.

\smallskip

We just need to verify the assumptions \ref{h6''} and \ref{h7''}.
Let us consider the new following assumptions: $\forall \eta>0,$
\begin{gather*}
\mathbb{P}\Big(\Big\{\sup_{\gep\leq t\leq T\wedge\tau_n^r}
\Big|B^{(n)}_i(t)-B^{(n)}_i(\gep)-\int_\gep^t
b_i(\X^{(n)}(s))\textnormal{d}
s\Big|\geq \eta\Big\}
\cap \Lambda_n(\alpha)\Big)
\underset{n\to\infty}{\longrightarrow}0
\tag*{(H6'')}\label{h6'} \\
%\forall \eta>0,\qquad
\mathbb{P}\Big(\Big\{
\sup_{\gep\leq t\leq T\wedge\tau_n^r}\Big|A_{ij}^{(n)}(t)-A_{ij}^{(n)}(\gep)
-\!\int_\gep^t\!\! a_{ij}(\X^{(n)}(s))%
\textnormal{d} s \Big|\geq\eta\Big\}\cap\Lambda_n(\alpha)\Big)
\underset{n\to\infty}{\longrightarrow}0,  \tag*{(H7'')}
\label{h7'}
\end{gather*}
where
\[
\Lambda_n(\alpha)=\Big(\Lambda_n^1(\alpha)\cup\Lambda_n^2(\alpha)\Big)^c.
\]
It is clear by Lemma \ref{restriction} that the assumptions \ref{h6'} and
\ref{h7'} yield the assumptions \ref{h6} and \ref{h7}.

\bigskip

\noindent\textbf{``Verification'' of assumption \ref{h6'}}

First, we change the supremum over real numbers into supremum over integers.
Cutting the integral over $b$ and using the fact that
on $\left[ \gep,(1-\varepsilon)\wedge\tau_{n}^{r}\right]$,
$\Xn_i(t)$ and $\Xn_i(t)-\Xn_j(t)$ are lower than $r$ and greater
than $\frac{1}{\sqrt{2n}}$, we obtain
\begin{align*}
\Big|\int_{\gep}^{t}&b_{i}(\Xn(s),s)\textnormal{d}s-\frac{1}{2n}%
\sum_{k=[2n\gep]}^{[2nt]-1}b_{i}(\Xn(\frac{k}{2n}),\frac{k}{2n})\Big|\\
&\leq\frac{1}{2n} \Big|\sum_{k=[2n\gep]}^{[2nt]-1}-\Xn_i(\frac{k}{2n})\int_{k}^{k+1}\frac
{1}{1-\frac{s}{2n}}\textnormal{d}s-\Xn_i(t)\int_{[2nt]}^{2nt}\frac
{1}{1-\frac{s}{2n}}\textnormal{d}s\\
&\hspace{1cm}-\sum_{k=[2n\gep]}^{[2nt]-1}\frac{-\Xn_i(\frac{k}{2n})}
{1-\frac{k}{2n}}\Big|
+\frac{(2p+1)r\sqrt{2n}}{2n}\\
&\leq
r\Big|\int_{2n\gep}^{2nt}\frac{1}{1-\frac{s}{2n}}\textnormal{d}s-\frac
{1}%
{2n}\sum_{k=[2n\gep]}^{[2nt]-1}\frac{1}{1-\frac{k}{2n}}\Big|
+{\scriptstyle\mathcal{O}}(1)\\
&\leq r\Big|\int_{\gep}^{t}\frac{1}{1-s}\textnormal{d}s-\frac{1}{2n}\sum
_{k=[2n\gep]}^{[2nt]-1}\frac{1}{1-\frac{k}{2n}}\Big|
+{\scriptstyle\mathcal{O}}(1)\\
& \leq{\scriptstyle\mathcal{O}}(1)
\end{align*}
where the
\textquotedblleft${\scriptstyle\mathcal{O}}(1)$\textquotedblright\
is uniform in $t$. As
$\Bn_i([t])=\Bn_i(t)+\mathcal{O}(\frac{1}{\sqrt{2n}})$, it
suffices to verify assumption \ref{h6'} to prove that for every $\eta>0$,
\begin{align*}
\mathbb{P}\Big(\Big\{
\sup_{\gep\leq k/2n\leq(1-\varepsilon)\wedge\tau_{n}^{r}}
\big|\Bn_i(\frac{k}{2n})-\Bn_i(\gep)-
\frac{1}
{2n}\sum_{j=[2n\gep]}^{k}b_{i}(&\Xn(\frac{j}{2n}),\frac{j}{2n})\big|\geq\eta\Big\}\\
&\hspace{.5cm}\cap\Lambda_n(\alpha)\Big)
\underset{n\rightarrow\infty}{\longrightarrow}0.
\end{align*}
For every function $h$, we denote by $\Delta_{k}h$ the increment of $h$
between the instants $\frac{k}{2n}$ and $\frac{k+1}{2n}$, {\it i.e.}
$$\Delta_{k}h=h(\frac{k+1}{2n})-h(\frac{k}{2n}).$$
With this notation, we have
\begin{align*}
\Delta_{k}\Bn_i
&=\mathbb{E}\left[\Delta_{k-1}\Xn_i\big|\Xn(\frac{k-1}{2n})\right]  \\
%&=\mathbb{P}\Big(  \Delta_{k-1}\Xn_i=1\big|\Xn(\frac{k-1}{2n})\Big)
%-\mathbb{P}\Big(  \Delta_{k-1}\Xn_i=-1\big|\Xn(\frac{k-1}{2n})\Big)\\
&=\sum_{\gep\in\{-1,1\}^p}
\frac{\gep_i}{\sqrt{2n}}\pr{\Delta_{k-1}\Xn_i=\frac{\gep}{\sqrt{2n}}
\big|\Xn(\frac{k-1}{2n})}
\end{align*}
and
\begin{align*}
\mathbb{P}\Big(\big(\Delta_{k}&\Xn_i\big)_{i}
=\big(\frac{\varepsilon_{i}}{\sqrt{2n}}\big)_{i}\big|
\big(\Xn_i(\frac{k}{2n})\big)_{i}=\big(x_{i}\big)_{i}\Big)\\
&=\mathbb{P}\left( \big(W^{(n)}_i(k+1)-
W_i^{(n)}(k)\big)_{i}=\big(\varepsilon_{i}\big)_{i}\big|
\big(X_i^{(n)}(\frac{k}{2n})\big)_{i}=\big(x_i\sqrt{2n}+1\big)_{i}\right)\\
&  =\frac{N(k,(x_i\sqrt{2n}-1)_{i}).N(2n-k-1,(x_i\sqrt{2n}-1+\varepsilon
_{i})_{i})}{N(k,(x_i\sqrt{2n}-1)_{i}).N(2n-k,(x_i\sqrt{2n}-1)_{i})}\\
&=\frac{N(2n-k-1,(x_i\sqrt{2n}-1+\varepsilon_{i})_{i})}{N(2n-k,(x_i\sqrt
{2n}-1)_{i})},
\end{align*}
where $N(m,e)$ denote the number of stars of length $m$ with wall condition
which end at $(e_1,\dots,e_p)$. Using (\ref{estim}),
\begin{align}
\mathbb{P}\Big(\big(\Delta_{k}&\Xn_i\big)_{i} =\big(\frac
{\varepsilon_{i}}{\sqrt{2n}}\big)_{i}\big|\big(\Xn_i%
(\frac{k}{2n})\big)_{i}=\big(x_i\big)_{i}\Big)\nonumber\\
%&=\prod\limits_{i=1}^{p}
%\frac{(\frac{1}{2}(2n-k+x_i\sqrt{2n}-1+p)!}{(\frac{1}%
%{2}(2n-k+x_i\sqrt{2n}-1+p+\frac{\varepsilon_{i}-1}{2})!}\nonumber\\
%&\phantom{=\prod\limits_{i=1}^{p}}\frac{(\frac{1}{2}(2n-k-x_i\sqrt{2n}-1+p-1)!}{(\frac{1}%
%{2}(2n-k-x_i\sqrt{2n}+1+p-1+\frac{\varepsilon_{i}-1}{2})!}
%\big(1+\frac{\varepsilon_{i}}{\sqrt
%{2n}x_i}\big)\frac{1}{2n-k+2i-2}\nonumber\\
%&\prod_{1\leq i<j\leq p}\big(1+\frac
%{\varepsilon_{j}-\varepsilon_{i}}{(x_j-x_i)\sqrt{2n}}\big)\big(1+\frac
%{\varepsilon_{j}+\varepsilon_{i}}{(x_j+x_i)\sqrt{2n}}\big)
%\nonumber\\
&  =\frac{1}{2^p}\prod\limits_{i=1}^{p}\big(1-\frac{\varepsilon_{i}x_i}%
{\sqrt{2n}(1-t)}+\frac{2p-1}{2n(1-t)}\big)\big(1+\frac
{\varepsilon_{i}}{x_i\sqrt{2n}+1}\big)\label{produit}\\
&  \phantom{=}\times\prod\limits_{1\leq i<j\leq p}\big(1+\frac
{\varepsilon_{j}-\varepsilon_{i}}{(x_j-x_i)\sqrt{2n}-2}\big)\big(1+\frac
{\varepsilon_{j}+\varepsilon_{i}}{(x_j+x_i)\sqrt{2n}+4}\big).\nonumber
\end{align}
We expand this product and we row these terms according to the number of
$1$, which appear in the product.

-- Choosing $1$ in every factor, we obtain of course $1$.

-- Choosing only one factor different from $1$, we obtain
\[\sum_{i=1}^{p}\big(\frac{-\varepsilon_{i}x_i%
}{\sqrt{2n}(1-t)}+\frac{\varepsilon_{i}}{x_i\sqrt{2n}}\big)+\sum_{1\leq
i<j\leq p}\big(\frac{\varepsilon_{j}-\varepsilon_{i}}{(x_j-x_i)\sqrt{2n}%
}+\frac{\varepsilon_{j}+\varepsilon_{i}}{(x_j+x_i)\sqrt{2n}}%
\big).
\]

-- Choosing two factors different from $1$, we obtain the sum of two
terms, one without epsilon and the other with the product of two different
epsilons:
\[
g(x_1,\dots,x_p)+\sum_{i\neq j}\gep_i\gep_jh_{i,j}(x_1,\dots,x_p)
\]
where $g$ and the $h_{i,j}$ are
$\mathcal{O}\big(\frac{1}{n^{2\alpha}}\big)$.

-- Choosing more than three factors different from $1$, we obtain some
terms which are $\mathcal{O}\big(\frac{1}{n^{3\alpha}}\big)$.

\smallskip

The expansion of the product (\ref{produit}) give
\begin{align}\label{transition}
\mathbb{P}\Big(\big(\Delta_{k}\Xn_i\big)_{i}=\big(&\frac
{\varepsilon_{i}}{\sqrt{2n}}\big)_{i}\big|\big(\Xn_i%
(\frac{k}{2n})\big)_{i}=\big(x_i\big)_{i}\Big)\nonumber\\
=&\frac{1}{2^{p}}\Big[1+\sum_{i=1}^{p}\big(\frac{-\varepsilon_{i}x_i%
}{\sqrt{2n}(1-t)}+\frac{\varepsilon_{i}}{x_i\sqrt{2n}}\big)\nonumber\\
&+\sum_{1\leq
i<j\leq p}\big(\frac{\varepsilon_{j}-\varepsilon_{i}}{(x_j-x_i)\sqrt{2n}%
}+\frac{\varepsilon_{j}+\varepsilon_{i}}{(x_j+x_i)\sqrt{2n}}%
\big)\\
&+g(x_1,\dots,x_p)+\sum_{i\neq j}\gep_i\gep_jh_{i,j}(x_1,\dots,x_p)
\Big]+\mathcal{O}\Big(\frac{1}{n^{3\alpha}}\Big).\nonumber
\end{align}
Multiplying the previous equality by $\gep_i/\sqrt{2n}$ and
summing over every $\gep\in\{-1,1\}^p$, all terms containing an odd number of
epsilon give $0$. Thus the increment of $B$ is
\begin{equation}\label{accroissement}
\Delta _{k}B^{(n)}_{i}=\frac{1}{2n}b_{i}\Big(\Xn(\frac{k}{2n}),\frac{k}{2n%
}\Big)+\mathcal{O}\Big(\frac{1}{n^{3\alpha+1/2}}\Big)
\end{equation}
Assumption \ref{h6'} is satisfied since
\begin{align*}
&\mathbb{P}\Big(\Big\{
\sup_{\gep\leq k/2n\leq (1-\varepsilon )\wedge \tau _{n}^{r}}\big|%
B_i^{(n)}(\frac{k}{2n})-B_i^{(n)}([2n\gep])\\
&\hspace{4.8cm}
-\frac{1}{2n}\sum_{j=[2n\gep]}^{k}b_{i}(
\Xn(\frac{j}{2n}),\frac{j}{2n})\big|%
\geq \eta\Big\}
\cap\Lambda_n(\alpha) \Big)\\
&=\mathbb{P}\Big(\Big\{
\sup_{\gep\leq k/2n\leq (1-\varepsilon )\wedge \tau _{n}^{r}}\big|%
\sum_{j=[2n\gep]}^{k}\Delta _{j}B^{(n)}_i-\frac{1}{2n}
b_{i}(\Xn(\frac{j}{2n}),\frac{j}{2n})%
\big|\geq \eta \Big\}\cap\Lambda_n(\alpha)\Big) \\
&\leq  \mathbb{P}\Big(\Big\{
\sum_{j=[2n\gep]}^{2n((1-\varepsilon )\wedge \tau _{n}^{r})}%
\big|\Delta _{j}B^{(n)}_i-\frac{1}{2n}
b_{i}(\Xn(\frac{j}{2n}),\frac{j}{2n})\big|%
\geq \eta \Big\}\cap\Lambda_n(\alpha)\Big) \\
&\leq\mathbb{P}\Big(\mathcal{O}(\frac{1}{n^{3\alpha-1/2}})\geq \eta
\Big\}\cap\Lambda_n(\alpha)\Big)%
\leq\mathbb{P}\Big(\mathcal{O}(\frac{1}{n^{3\alpha-1/2}})\geq \eta\Big)
\underset{n\rightarrow \infty}{\longrightarrow }0
\end{align*}%
and $\alpha>1/6$.

\bigskip

\noindent\textbf{``Verification'' of assumption \ref{h7'}}

First, we have
\begin{align*}
M_i^{(n)}(\frac{k+1}{2n}&)M_j^{(n)}(\frac{k+1}{2n})\\
%=&\big(\X_i^{(n)}(k)-B^{(n)}_i(k)-\Delta_k\X_i^{(n)}+\Delta_kB^{(n)}_i\big)\\
%&\big(\X_j^{(n)}(k)-B_j^{(n)}(k)-\Delta_k\X_j^{(n)}+\Delta_kB_j^{(n)}\big)\\
=&M_i^{(n)}(\frac{k}{2n})M_j^{(n)}(\frac{k}{2n})
+\big(\X_i^{(n)}(\frac{k}{2n})-B^{(n)}_i(\frac{k}{2n})\big)\big(
\Delta_kB_j^{(n)}-
\Delta_k\X_j^{(n)}\big)\\
&+\big(\X_j^{(n)}(\frac{k}{2n})-B_j^{(n)}(\frac{k}{2n})\big)\big(\Delta_kB^{(n)}_i-\Delta_k\X_
i^{(n)}\big)\\
&+\big(\Delta_kB^{(n)}_i-\Delta_k\X_i^{(n)}\big)
\big(\Delta_kB_j^{(n)}-\Delta_k\X_j^{(n)}\big)
\end{align*}
and by definition of $B$,
$$\esp{\Delta_kB^{(n)}-\Delta_k\Xn\big|\mathcal{F}_k}=0,$$
hence
\begin{align*}
\esp{\Delta_kM_i^{(n)}M_j^{(n)}\big|\mathcal{F}_k}
%&=\esp{\big(\Delta_kB^{(n)}_i-\Delta_k\X_i^{(n)}\big)
%\big(\Delta_kB_j^{(n)}-\Delta_k\Xn_j\big)\big|\mathcal{F}_k}\\
%&
=\esp{\Delta_k\Xn_i\Delta_k\Xn_j\big|\mathcal{F}_k}-
\Delta_kB^{(n)}_i\Delta_kB_j^{(n)}.
\end{align*}
On $\Lambda_n(\alpha)$ and for
$t\in[\gep,((1-\gep)\wedge\tau_n^r)]$,
$X_i^{(n)}(t)$ is lower than $r$
and greater than $n^{\alpha-1/2}$ thus
\[
b_{i}(\Xn(\frac{k}{2n}),\frac{k}{2n})=\mathcal{O}(n^{1/5}).\]
Using (\ref{accroissement}),
$$\sum_{l=[2n\gep]}^{k-1}\Delta _{l}B^{(n)}_i\Delta _{l}B_j^{(n)}
=\mathcal{O}(\frac{1}{n^{1-2\alpha}}).$$
This estimate yields
\begin{align*}
A^{(n)}_{i,j}(\frac{k}{2n})-A^{(n)}_{i,j}(\gep)&=\sum_{l=[2n\gep]}^{k-1}
\esp{\Delta_lM_i^{(n)}M_j^{(n)}\big|\mathcal{F}_l}\\
%&=\sum_{l=[2n\gep]}^{k-1}\esp{\Delta_lX_i^{(n)}\Delta_lX_j^{(n)}\big|
%\mathcal{F}_l}
%-\sum_{l=[2n\gep]}^{k-1}\Delta_lB^{(n)}_i\Delta_lB_j^{(n)}\\
&=\sum_{l=[2n\gep]}^{k-1}\esp{\Delta_lX_i^{(n)}\Delta_lX_j^{(n)}\big|
\mathcal{F}_l}
+\mathcal{O}\Big(\frac1{n^{1-2\alpha}}\Big).
\end{align*}
%It remains to estimate the last sum.
We consider two cases whether
$i$ equals to $j$ or not.

If $i=j$ then, since $(\Delta_lX_i^{(n)})^2=1/2n$, we have
\begin{align*}
A^{(n)}_{i,j}(\frac{k}{2n})-A^{(n)}_{i,j}(\gep)
&=\frac{k-[2n\gep]}{2n} +\mathcal{O}\Big(\frac1{n^{1-2\alpha}}
\Big)\\
&=\int_\gep^{k/2n}a_{i,i}(t)\dd t
+\mathcal{O}\Big(\frac1{n^{1-2\alpha}}\Big).
\end{align*}

If $i<j$ then
$$\esp{\Delta_lX_i^{(n)}\Delta_lX_j^{(n)}\big|\mathcal{F}_l}
=\mathcal{O}\Big(\frac1{n^2}\Big).$$
Indeed, if we denote  by
$$p_k(\gep_i,\gep_j,x)
=\pr{\Delta_lX_i^{(n)}=\gep_i,\Delta_lX_j^{(n)}=\gep_j\big|X_n(\frac{l}{2n})=x},$$
for $\gep_i$ and $\gep_j$ in $\{-1,1\}$ and for
$x=(x_i)$ a vector of $\mathbb{R}^p$, then %alors
$$\esp{\Delta_l X_i^{(n)}\Delta_lX_j^{(n)}\big|X_n(\frac{l}{2n})=x}
=\sum_{\gep_i,\gep_j\in\{-\frac{1}{\sqrt{2n}},\frac{1}{\sqrt{2n}}\}}
\gep_i\gep_jp_l(\gep_i,\gep_j,x).$$
Now, using (\ref{transition}) and the order of magnitude of the functions
$g$ and $h$, we have the following estimate
\begin{align*}
p_l(\gep_i,\gep_j,x)
&=\sum_{\substack{k\in\{1,\dots,p\},k\neq i,j\\
\gep_k\in\{-\frac{1}{\sqrt{2n}},\frac{1}{\sqrt{2n}}\}}}
\pr{\Delta_lX^{(n)}=\gep\big|X_n(\frac{k}{2n})=x }\\
&=\frac{1}{2^p}\sum_{\substack{k\in\{1,\dots,p\},k\neq i,j\\
\gep_k\in\{-\frac{1}{\sqrt{2n}},\frac{1}{\sqrt{2n}}\}}}
\Big[1+\sum_{h=1}^p\big(\frac{-\gep_hx_h}{\sqrt{2n}(1-t)}
+\frac{\gep_h}{x_h\sqrt{2n}}\big)\\
&%\hspace{3.48cm}
+\sum_{1\leq r<s\leq p}\big(
\frac{\gep_r-\gep_s}{(x_r-x_s)\sqrt{2n}}
+\frac{\gep_r+\gep_s}{(x_r+x_s)\sqrt{2n}}\big)
\Big]+\mathcal{O}\Big(\frac1{n^{2\alpha}}\Big)\\
&=\frac14+\frac1{4\sqrt{2n}}\Big[\frac{-\gep_ix_i}{1-l/n}+\frac{\gep_i}{x_i}
+\frac{-\gep_jx_j}{1-l/n}+\frac{\gep_j}{x_j}+\frac{\gep_i-\gep_j}{x_i-x_j}
+\frac{\gep_i+\gep_j}{x_i+x_j}\\
&%\hspace{2.5cm}
+\sum_{\substack{s\in\{1,\dots,p\}\\s\neq i,j}}\big(
\frac{\gep_i}{x_i-x_s}+\frac{\gep_i}{x_i+x_s}+\frac{\gep_j}{x_j-x_s}
+\frac{\gep_j}{x_j+x_s}\big)\Big]+\mathcal{O}\Big(\frac{1}{n^{2\alpha}}\Big).
\end{align*}
A simple computation leads to
$$\esp{\Delta_lX_i^{(n)}\Delta_lX_j^{(n)}\big|\Xn(\frac{l}{2n})}
=\sum_{\gep_i,\gep_j\in\{-\frac{1}{\sqrt{2n}},\frac{1}{\sqrt{2n}}\}}
\gep_i\gep_jp_l(\gep_i,\gep_j,x)
=\mathcal{O}\Big(\frac1{n^{1+2\alpha}}\Big).$$
Thanks to this estimate, we can write, for $i\neq j$,
\begin{align*}
A_{i,j}^{(n)}(\frac{k}{2n})-A_{i,j}^{(n)}(\gep)
&=\mathcal{O}\Big(\frac{1}{n^{2\alpha}}\Big)
+\mathcal{O}\Big(\frac{1}{n^{1-2\alpha}}\Big)\\
&=\int_\gep^{k/2n}a_{i,j}(t)\dd t+\mathcal{O}\Big(\frac{1}{n^{2\alpha}}\Big).
\end{align*}

In both cases, we have shown
$$A_{i,j}^{(n)}(\frac{k}{2n})-A^{(n)}_{i,j}(\gep)-\int_\gep^{k/2n}a_{i,j}(t)\dd
t=\mathcal{O}\Big(\frac{1}{n^{2\alpha}}\Big).$$
This equality completes the checking of Assumption \ref{h7'} and
the proof of Theorem \ref{convergence} as well.
\end{proof}

\section{Watermelons without wall condition: proofs}
\label{preuvesm}

\subsection{Properties}\label{propsm}

\begin{prop}\label{bessel sm}
Let $\widehat{X}=(\widehat{X}_1,\dots,\widehat{X}_p)$ be a solution to the SDE \ref{EDSsm}. The
Euclidean norm of
   $\widehat{X}$ is a Bessel bridge
with dimension $p^2$.
\end{prop}

\begin{proof}
Using the It\^o's formula, we prove that $U=\|\widehat{X}\|^2$ satisfies the following
equation:
$$\dd U=\Big(\frac{-2U}{1-t}+p^2\Big)\dd t+\sqrt{U}\dd B.$$
This equation is the SDE from a square Bessel bridge with dimension $p^2$.
\end{proof}

\begin{prop}\label{vander sm}
If $\widehat{X}=(\widehat{X}_1,\dots,\widehat{X}_p)$ is a solution of the SDE \ref{EDSsm}, then we have
$$
\pr{\forall t \in \left]0;1\right[,\
    -\infty<\widehat{X}_1(t)<\cdots<\widehat{X}_p(t)<+\infty}=1.
$$
\end{prop}

\begin{proof}
Let $\widehat{X}$ be a solution of SDE \ref{EDSsm}.
By Proposition \ref{bessel sm}, the norm of $\widehat{X}$ is a Bessel bridge
   with dimension $p^2$, since $\widehat{X}$ is finite on $[0,1]$.
It remains to show that the branches of $\widehat{X}$ do not touch each
other. For this we use the same method as in the previous section.
Set $\gep >0$, $F(x_1,\dots,x_p)=\ln\big(\prod_{i<j}(x_j-x_i)\big)$ and
$U(t)=F(\widehat{X}(t))$ for $t\in[\gep,(1-\gep)]$. The partial derivatives of $F$
are given by
\[
\partial_iF(x)=\sum_{j\neq i}\frac{1}{x_i-x_j}
\qquad\text{ and }\qquad
\partial_i^2F(x)=\sum_{j\neq i}\frac{-1}{(x_i-x_j)^2}.
\]
By Girsanov Theorem, it exists a probability measure $\mathbb{Q}$ under
which $U$ satisfies 
\[
\dd U=\sum_{i=1}^p\big((\partial_iF(\widehat{X}))^2+\frac12\partial_i^2F(\widehat{X})\big)\dd s
+\sum_{i=1}^p\partial_iF(\widehat{X})\dd B_i.
\]
Now
\[
\sum_{i=1}^p(\partial_iF(x))^2+\frac12\partial_i^2F(x)
=\frac12\sum_{i\neq j}\frac{1}{(x_j-x_i)^2}
+\sum_{i\neq j\neq k}\frac{1}{(x_i-x_j)(x_i-x_k)}=0
\]
thus $U$ is a time-changed Brownian motion, $\mathbb{Q}$-almost surely
finite on $[\gep,(1-\gep)]$  and thus $\mathbb{P}$-almost
surely for every $\gep>0$. This completes the proof of proposition
\ref{vander sm}.
\end{proof}

\begin{proof}[Proof of Theorem \ref{th a p branches sm}]
We use the same method as in the proof of Theorem \ref{th a p
    branches}. Since the coefficients of the SDE \ref{EDSsm} are locally
Lipschitz, we have the uniqueness of the solution up to its explosion
time. Proposition \ref{vander sm} implies that this explosion time equals
to $1$. Thus we have  the uniqueness of \ref{EDSsm} over $[\gep,1-\gep]$
for every $\gep>0$ and we proceed the same way to obtain the uniqueness
over $[0,1]$.

The existence of a solution of \ref{EDSsm} results from the fact that the
renormalized sequence of the $(p,2n)$-watermelons has a subsequence
which converges
to a solution of \ref{EDSsm}. This point will be proved during the proof of
Theorem \ref{convergence sm}.
\end{proof}

\subsection{1-dimensional distribution}

\begin{prop}\label{densite sm}
Let $\widehat{X}^{(n)}$ be a renormalized $(p,2n)$-watermelon without
wall condition an let $t$ be in $[0;1]$. $\widehat{X}^{(n)}(t)$ converges in
distribution to $\sqrt{2t(1-t)}\ \widehat{\Lambda}$.

\end{prop}

\begin{proof}
We denote by $\widehat{f}(t;x)$ the density function of 
 $\sqrt{2t(1-t)}\ \widehat{\Lambda}$ {\it i.e.}
\begin{equation*}
\widehat{f}(t;x)=\frac{2^{-p/2}}{\pi^{p/2}\big(t(1-t)\big)^{p^2/2}
\prod_{i=1}^{p-1}i!}
\prod_{1\leq i<j\leq p}(x_j-x_i)^2
\ \mathrm{e}^{-\frac{\|x\|^2}{2t(1-t)}}1%
\hspace{-0.9mm}\mathrm{l}_{x_1\leq \cdots\leq x_p}.
\end{equation*}
Let $u$ and $v$ be two vectors in $\mathbb{R}^p$ such that $u<v$ and let
$t$ be in $\left[0;1\right]$.
To prove Proposition \ref{densite sm}, it suffices to show that 
\begin{equation*}
\mathbb{P}\left(\widehat{X}^{(n)}(t)\in\left[u,v\right]\right)
\underset{n\to\infty}{\longrightarrow}\int_{[u;v]}\widehat{f}(t;x)
\mathnormal{d}x.
\end{equation*}
We denote by $\widehat{N}(m,e)$ the number of star without wall condition of length m
with $p$ branches ending at $(e_1,\dots,e_p)$. By Theorem 1. in
\cite{viennot}, we know that
\begin{align}\label{estimsm}
\widehat{N}(m,e)=2^{\binom{p}{2}}
\prod_{i=1}^p\frac{(m-i+p)!}
{\big(\pfrac{12}(m+e_i)\big)!\big(\pfrac{12}(m-e_i)+p-1\big)!}
\prod_{1\leq i<j\leq p}(e_j-e_i).
\end{align}
Let $u,v\in\mathbb{R}^p$ and $m=[2nt]$ where $t\in\left[0;1\right]$, we have
$$
\pr{\widehat{X}^{(n)}(t)\in\left[u,v\right]}=
\sum_{\substack{u\leq x\leq v\\x\sqrt{2n}\in\mathbb{Z}}}
\pr{\widehat{W}^{(n)}([2nt])=x\sqrt{2n}}.
$$
Cutting the watermelons in two stars in the same way as in the proof
of Proposition \ref{densite sm}, we can write that
\begin{align*}
\pr{W^{(n)}([2nt])=x\sqrt{2n}}
&=\frac{\widehat{N}(m,x\sqrt{2n})\widehat{N}(2n-m,x\sqrt{2n})}{\widehat{N}(2n,2i-2)}\\
&=\widehat{c}_{p,n}\prod_{i=1}^p\frac{\widehat{D}(m,x_i,p)\widehat{D}(2n-m,x_i,p)}
{\widehat{D}(2n,\frac{2i-2}{\sqrt{2n}},p)}
\end{align*}
where
\[
\widehat{c}_{p,n}=\frac{2^{-\binom{p}{2}}
\Big[\prod_{1\leq i<j\leq p}(x_j-x_i)\sqrt{2n}\Big]^2}
{\prod_{1\leq i<j\leq p}2(j-i)}
=\frac{n^{\binom{p}{2}}
\Big[\prod_{1\leq i<j\leq p}(x_j-x_i)\Big]^2}
{2^{\binom{p}{2}}\prod_{i=1}^{p-1}i!}
\]
and
\[
\widehat{D}(m,x_i,p)=\frac{(m-i+p)!}
{(\pfrac{12}(m+x_i\sqrt{2n}))!(\pfrac{12}(m-x_i\sqrt{2n})+p-1)!}\cdot
\]
Lemma \ref{lemme} Yields the following estimate
\[
\widehat{D}(m,x_i,p)=\frac1{\sqrt{\pi}}(nt)^{-i+\pfrac{12}}2^{m+p-i}
\,\textnormal{e}^{-{x_i}^2/(2t)}
\Big\{1+\mathcal{O}\Big(\frac{1}{\sqrt{n}}\Big)\Big\}
\]
and thus
\begin{align*}
&\frac{\widehat{D}(m,x_i,p)\widehat{D}(2n-m,x_i,p)}{\widehat{D}(2n,\frac{2i-1}{\sqrt{2n}},i-1+p)}
=\frac{2^{p-i}}{\sqrt{\pi}}\big(nt(1-t)\big)^{\pfrac{12}-i}
\ \textnormal{e}^{\frac{-{x_i}^2}{2t(1-t)}}
\Big\{1+\mathcal{O}\Big(\frac{1}{\sqrt{n}}\Big)\Big\}
\end{align*}
and
\begin{multline*}
\pr{\widehat{W}^{(n)}([2nt])=x\sqrt{2n}}\\
=\Big(\frac{2}{n}\Big)^{p/2}\frac{2^{-p/2}}{\pi^{p/2}\prod_{i=1}^{p-1}i!}
\frac{\prod_{1\leq i<j\leq p}(x_j-x_i)^2}{\big(t(1-t)\big)^{p^2/2}}
\textnormal{e}^{-\frac{\|x\|^2}{2t(1-t)}}
\Big\{1+\mathcal{O}\Big(\frac{1}{\sqrt{n}}\Big)\Big\}.
\end{multline*}
Finally, we obtain
\begin{multline*}
\pr{\widehat{X}^{(n)}(t)\in\left[u,v\right]}\\
=\Big(\frac{2}{n}\Big)^{p/2}
\sum_{\substack{u\leq x\leq v\\x\sqrt{2n}\in\mathbb{N}}}
\frac{2^{-p/2}}{\pi^{p/2}\prod_{i=1}^{p-1}i!}
\frac{\prod_{1\leq i<j\leq p}(x_j-x_i)^2}{\big(t(1-t)\big)^{p^2/2}}
\textnormal{e}^{-\frac{\|x\|^2}{2t(1-t)}}
\Big\{1+\mathcal{O}\Big(\frac{1}{\sqrt{n}}\Big)\Big\}.
\end{multline*}
As the $\mathcal{O}$ are uniform in $x$, the previous sum is a Riemann sum
which converges to $\int_{[u;v]}\widehat{f}(t;x)\dd x$.
\end{proof}

\subsection{Proof of Theorem \ref{convergence sm}}
Let $\big(\widetilde{X}^{(n)}(t)\big)_{t\in[0;2n]}$ be the piecewise continuous function such that for every
$k\in\oc0;2n\fc$ and every $i\in\oc1;p\fc$,
\[
\X^{(n)}_i(k)=
\widehat{X}^{(n)}_i(k)+\frac{1}{\sqrt{2n}}
\]
As in the proof of Theorem \ref{convergence}, it is enough to prove the
convergence of the process $\X^{(n)}$ to the process $\widehat{X}$, solution of the
SDE \ref{EDSsm} over $[\gep,1-\gep]$ for every $\gep\in]0,1/2[$.
For this, we shall use Theorem \ref{ethier}.

\medskip

Let the processes $A$, $B$ and $M$, the stopping times $T$ and $\tau_n^r$
and the functions $\sigma$ and $a$ de defined as in the proof of
Theorem \ref{convergence}. We just replace $f$ by $\widehat{f}$
in the definition of property $\mathcal{P}$ and defined the function
$b$ by
\[
b(t;x)=(b_i(t,x))_{1\leq i\leq p}
=\Big(\frac{-x_i}{1-t}+\sum_{j=1}^p\frac{1}{x_i-x_j}\Big)_{1\leq i\leq p}.
\]
Almost all assumptions of Theorem \ref{ethier} are easily verified.
Indeed, it is
enough to check the assumption \ref{h6'} and \ref{h7'}.
First, let us give a lemma similar to Lemma \ref{restriction}.

For $\gep\in]0;1[$ and $\alpha\in]0;1/4[$, we define the set
\begin{multline*}
\Lambda_n(\alpha)=\Big\{\exists t\in [\gep, ((1-\gep)\wedge\tau_n^r)]
,\ \exists i<j\in\oc1,p\fc \\
\text{ such that } \widehat{W}_j^{(n)}([2nt])-\widehat{W}_i^{(n)}([2nt])\leq n^\alpha
\Big\}.
\end{multline*}
We have
\begin{lem}\label{restriction sm}
\[\pr{\Lambda_n(\alpha)}\underset{n\to\infty}{\To}0.\]
\end{lem}

\begin{proof}
Proposition \ref{densite sm} yields the following equality
\[
\pr{\X^{(n)}(t)\in\left[u,v\right]}
=\frac{1}{(2n)^{p/2}}
\sum_{\substack{u\leq x\leq v\\
x\sqrt{2n}\in\mathbb{Z}}}
\widehat{f}(t;x)\Big\{1+\mathcal{O}\Big(\frac{1}{\sqrt{n}}\Big)\Big\}.
\]
Setting an upper bound to $\widehat{f}(t;x)$, we obtain
\[
\pr{\widehat{W}^{(n)}_j([2nt])-\widehat{W}_i^{(n)}([2nt])\in[0,n^\alpha]}
=\mathcal{O}\big(n^{3(\alpha-1/2)}\big)
\]
and so
\[
\pr{\Lambda_n(\alpha)}=\mathcal{O}\big(n^{2\alpha-1/2}\big).
\]
This estimate completes the proof since $\alpha<1/4$.
\end{proof}

\bigskip

\noindent {\bf ``Verification'' of assumption \ref{h6'}}

First, we change the supremum over real numbers into a supremum over
integers. The assumption \ref{h6'} becomes
\begin{multline*}
\mathbb{P}\Big(\Big\{
\sup_{\gep\leq
k/2n\leq(1-\varepsilon)\wedge\tau_{n}^{r}}\big|\Bn_i(\frac{k}{2n})-
\Bn_i(\gep)-
\frac{1}{2n}\sum_{j=1}^{k}b_{i}(\Xn(\frac{j}{2n}),\frac{j}{2n})\big|
\geq\eta\Big\}\\
\cap\Lambda_n(\alpha)^c\Big)
\underset{n\rightarrow\infty}{\longrightarrow}0.
\end{multline*}
Let us recall that for a function $h$,
$\Delta_kh=h(\frac{k+1}{2n})-h(\frac{k}{2n})$. By definition
of $B$, we can write
\[
\Delta_kB_i^{(n)}=\sum_{\gep\in\{-1,1\}}\frac{\gep_i}{\sqrt{2n}}
\pr{\Delta_{k-1}\Xn_i=\frac{\gep_i}{\sqrt{2n}}\big|\Xn(\frac{k-1}{2n})}
\]
and for $\gep\in\{-1,1\}^p$ and $x\in\{x_1<\cdots<x_p\in\mathbb{R}\}$
\[
\mathbb{P}\Big(\Delta_{k-1}\X^{(n)}=\frac{\gep}{\sqrt{2n}}\big|
\Xn(\frac{k-1}{2n})=x\Big)
=\frac{\widehat{N}(2n-k-1,(x_{i}\sqrt{2n}-1+\varepsilon_{i})_{i})}{\widehat{N}(2n-k,(x_{i}\sqrt
{2n}-1)_{i})}.
\]
where $\widehat{N}$ denotes the number of stars without wall condition. Using
(\ref{estimsm}) and the fact that on $\Lambda_n(\alpha)$, $x>n^{\alpha-1/2}$
and that $k\in\oc 2n\gep,2n((1-\gep)\wedge\tau_n^r)\fc$, we obtain the
following estimate
\begin{align*}
\mathbb{P}\Big(&\big(\Delta_{k}\Xn_i\big)_{i} =\big(\frac
{\varepsilon_{i}}{\sqrt{2n}}\big)_{i}\big|\big(\Xn_i%
(\frac{k}{2n})\big)_{i}=\big(x_{i}\big)_{i}\Big)\\
=&\frac{1}{2^p}\prod\limits_{i=1}^{p}\big(1-\frac{\varepsilon_{i}x_i}%
{\sqrt{2n}(1-t)}+\frac{2p(1+\gep_i)-1}{2n(1-t)}\big)
\prod\limits_{1\leq i<j\leq p}
\big(1+\frac{\varepsilon_{j}-\varepsilon_{i}}{(x_{j}-x_{i})\sqrt{2n}-2}\big)\\
=&\frac{1}{2^{p}}\Big[
1
+\sum_{i=1}^{p}\frac{-\varepsilon_{i}x_{i}}{\sqrt{2n}(1-t)}
+\sum_{1\leq i<j\leq p}
\frac{\varepsilon_{j}-\varepsilon_{i}}{(x_{j}-x_{i})\sqrt{2n}}
+g(x)+\sum_{i\neq j}\gep_i\gep_jh_{i,j}(x)\Big]\\
&\hspace{9.8cm}
+\mathcal{O}\Big(\frac{1}{n^{3\alpha}}\Big)
\end{align*}
where $g$ and the $h_{i,j}$ are $\mathcal{O}\big(n^{-2\alpha}\big)$.
Hence, we have
\begin{align}\label{accroissement sm}
\Delta_kB^{(n)}_i=\frac{1}{2n}b_i\Big(\Xn(\frac{k}{2n}),\frac{k}{2n}\Big)
+\mathcal{O}\Big(\frac{1}{n^{3\alpha+1/2}}\Big)
\end{align}
and so
%\begin{align*}
\begin{multline*}
\mathbb{P}\Big(\Big\{
\sup_{k/2n\leq(1-\varepsilon)\wedge\tau_{n}^{r}}
\big|\Bn_i(\frac{k}{2n})-\Bn_i(\gep)
-\frac{1}{2n}\sum_{j=[2n\gep]}^{k}b_{i}(\Xn(\frac{j}{2n}),\frac{j}{2n})
\big|\geq\gep\Big\}\\
\cap\Lambda_n(\alpha)^c\Big)
\leq
\mathbb{P}\Big(\mathcal{O}\big(\frac{1}{n^{3\alpha+1/2}}\big)\geq\eta\Big)
\underset{n\rightarrow\infty}{\longrightarrow}0
\end{multline*}
%\end{align*}
since $\alpha$ is greater than $1/6$.

\bigskip

\noindent{\bf ``Verification'' of assumption \ref{h7'}}

The definitions being the same as in Theorem \ref{convergence}, we
still have
\[
\esp{\Delta_kM_i^{(n)}M_j^{(n)}\big|\mathcal{F}_k}
=\esp{\Delta_k\Xn_i\Delta_k\Xn_j\big|\mathcal{F}_k}-
\Delta_kB^{(n)}_i\Delta_kB_j^{(n)}.
\]
Recall that on $\Lambda_n(\alpha)$ and for
$k\in\oc 2n\gep,((1-\gep)\wedge\tau_n^r)2n\fc$,
we have
$$b_{i}(\widetilde{X}(\frac{k}{2n}),\frac{k}{2n})=\mathcal{O}(n^{1/5})$$
and therefore, by (\ref{accroissement sm}),
$$\sum_{l=[2n\gep]}^{k-1}\Delta _{l}B^{(n)}_i\Delta _{l}B_j^{(n)}
=\mathcal{O}(\frac{1}{n^{1-2\alpha}}).$$
Hence, we have

\begin{align*}
A^{(n)}_{i,j}(k)
=\sum_{l=[2n\gep]}^{k-1}\esp{\Delta_l\Xn_i\Delta_l\Xn_j\big|
\mathcal{F}_l}
+\mathcal{O}\Big(\frac1{n^{1-2\alpha}}\Big).
\end{align*}

\medskip

If $i=j$, then, since $(\Delta_l\Xn_i)^2=1/2n$, we have
$$A^{(n)}_{i,j}(\frac{k}{2n})=\frac{k-[2n\gep]}{n} +\mathcal{O}\Big(\frac1{n^{1-2\alpha}}
\Big)
=\int_\gep^{k/2n}a_{i,i}(t)\dd t +\mathcal{O}\Big(\frac1{n^{1-2\alpha}}\Big).$$

\medskip

If $i<j$, then the equality
$$\esp{\Delta_l\Xn_i\Delta_l\Xn_j\big|\mathcal{F}_l}
=\mathcal{O}\Big(\frac1{n^{1+2\alpha}}\Big)$$
yields
$$A^{(n)}_{i,j}(k)-A^{(n)}_{i,j}(\gep)
=\int_\gep^{k/2n}a_{i,j}(t)\dd
t+\mathcal{O}\Big(\frac{1}{n^{2\alpha}}\Big).
$$

\medskip

In both cases, we have
$$A^{(n)}_{i,j}(\frac{k}{2n})-A^{(n)}_{i,j}(\gep)-\int_\gep^{k/2n}a_{i,j}(t)\dd t
=\mathcal{O}\Big(\frac{1}{n^{2\alpha}}\Big).$$
This complete the verification of the assumption \ref{h7'} and
the proof of Theorem \ref{convergence sm} too.

%%%%%%%%%%%%%%%%%%%%%%%%%%%%%%%%%%%%%%%%%%%%%%%%%%%%%%%%%%%%%%%%%%%%%%%%%%
%%%%%%%%%%%%%%%%%%%%%%%%%%%%%%%%%%%%%%%%%%%%%%%%%%%%%%%%%%%%%%%%%%%%%%%%%%

\section{Some moments of asymptotic watermelons}

In this section, we give in both first propositions the moments at
time $t$ of the continuous $2$-watermelons with and
without wall condition. We give in the sequel for a continuous
$p$-watermelon $X$ the moments of the elementary symmetric polynomials
of $X_1,\dots,X_p$.

\begin{prop}\label{moment avec mur}
Let $(X_1,X_2)$ be a continuous $2$-watermelon with
wall condition. We have for every integer $k\in\mathbb{N}^*$

\begin{align*}
&\frac{3\pi\esp{X_2(t)^{2k}}}{\big(t(1-t)\big)^k}=
%\\&
2(3-k)(k+1)!+\frac{(k^2+k+3)(2k+2)!}{(k+1)!2^{k+1}}
\Big[\pi-2\sum_{j=1}^{k+1}
\frac{2^{j}}{j\binom{j}{2j}}\Big]
%\frac{\big(t(1-t)\big)^k}{3\pi}
,\\
&\frac{3\pi\esp{X_1(t)^{2k}}}{\big(t(1-t)\big)^k}=
-2(3-k)(k+1)!+\frac{(k^2+k+3)(2k+2)!}{(k+1)!2^{k+1}}
\Big[\pi+2\sum_{j=1}^{k+1}
\frac{2^{j}}{j\binom{j}{2j}}\Big],\\
&\frac{3\sqrt{\pi}\esp{X_1(t)^{2k-1}}}{\big(t(1-t)\big)^{(2k-1)/2}}=
\frac{(2k+2)!(14-4k)}{2^{2k+3}(k+1)!}+(4k^2+11)2^{k-1}k!
\big(\sqrt{2}-\sum_{j=0}^k\frac{\binom{j}{2j}}{2^{3j}}\big),\\
&\frac{3\sqrt{\pi}\esp{X_2(t)^{2k-1}}}{\big(t(1-t)\big)^{(2k-1)/2}}=
-\frac{(2k+2)!(14-4k)}{2^{2k+3}(k+1)!}+(4k^2+11)2^{k-1}k!
\sum_{j=0}^k\frac{\binom{j}{2j}}{2^{3j}}.
\end{align*}

\end{prop}

\begin{prop}\label{moment sans mur}
Let $(\widehat{X}_1,\widehat{X}_2)$ be a continuous $2$-watermelon
without wall condition. We have for every integer
$k\in\mathbb{N}$
\begin{align*}
\esp{\widehat{X}_2(t)^{2k}}&=\esp{\widehat{X}_1(t)^{2k}}=
\frac{(2k)!(k+1)}{2^{k}k!}\big(t(1-t)\big)^{k},\\
\esp{\widehat{X}_2(t)^{2k+1}}&=-\esp{\widehat{X}_1(t)^{2k+1}}\\
&=\Big[\frac{(2k+2)!}{2^{2k+1}(k+1)!}
+(2k+3)2^kk!\sum_{j=0}^k
\frac{\binom{j}{2j}}{2^{3j}}\Big]
\frac{\big(t(1-t)\big)^{(2k+1)/2}}{2\sqrt{\pi}}.
\end{align*}
\end{prop}

\noindent{\bf Remark :}
The following table give the values of the first moments of
continuous $2$-watermelons:

\[
\begin{array}{|c|@{\,}|c|c|@{\,}|c|c|}
\hline
&\multicolumn{2}{c|@{\,}|}{\vbox{\hsize=5cm
\begin{tabular}{c} $2$-watermelon \\ without wall condition\end{tabular}}}
&\multicolumn{2}{c|}{\quad\vbox{\hsize=5cm
\begin{tabular}{c}$2$-watermelon\\with wall condition\end{tabular}}}\\
\hline
 k
&\vbox{\hsize=2.4cm $\frac{2\pi\esp{\widehat{X}_1(t)^k}}{(t(1-t))^{k/2}} $}
&\vbox{\hsize=2.4cm $\frac{2\pi\esp{\widehat{X}_2(t)^k}}{(t(1-t))^{k/2}} $}
&\vbox{\hsize=2.4cm $\frac{3\pi\esp{X_1(t)^k}}{(t(1-t))^{k/2}} $}
&\vbox{\hsize=2.4cm $\frac{3\pi\esp{X_2(t)^k}}{(t(1-t))^{k/2}} $}\\
\hline\hline
\scriptstyle1&\scriptstyle-4\sqrt{\pi}&\scriptstyle4\sqrt{\pi}&\scriptstyle(15\sqrt{2}-15)\sqrt{\pi}&\scriptstyle15\sqrt{\pi}\\
\hline
\scriptstyle2&\scriptstyle4\pi&\scriptstyle4\pi&\scriptstyle15\pi-32&
\scriptstyle15\pi+32\\
\hline
\scriptstyle3&\scriptstyle-14\sqrt{\pi}&\scriptstyle14\sqrt{\pi}
&\scriptstyle(108\sqrt{2}-\frac{279}{2})\sqrt{\pi}&\scriptstyle\frac{279}{2}\sqrt{\pi}\\
\hline
\scriptstyle4&\scriptstyle18\pi&\scriptstyle18\pi&\scriptstyle135\pi-384&\scriptstyle135\pi+384\\
\hline
\scriptstyle5&\scriptstyle-79\sqrt{\pi}&\scriptstyle79\sqrt{\pi}
&\scriptstyle(1128\sqrt{2}-\frac{6213}{4})\sqrt{\pi}&\scriptstyle\frac{6213}{4}\sqrt{\pi}\\
\hline
\scriptstyle6&\scriptstyle120\pi&\scriptstyle120\pi&\scriptstyle1575\pi-4800&\scriptstyle1575\pi+4800\\
\hline
\end{array}
\]

\begin{prop}\label{sym1}
Let $(X_i)_{1\leq i\leq p}$ be a continuous $p$-watermelon with
wall condition and $\Sigma_{k,p}$ be the $k$-th elementary symmetric
polynomial of $({X_i}^2)_{1\leq i\leq p}$,
{\it i.e.}
$$\Sigma_{k,p}=\sum_{1\leq i_1<\dots< i_k\leq p}
{X_{i_1}}^2\dots {X_{i_k}}^2.$$
We have
\[
\esp{\Sigma_{2k,p}}=\frac{(2p+1)!}{(2p+1-2k)!2^kk!}\big(t(1-t)\big)^k.
\]
\end{prop}
\begin{prop}\label{sym2}
Let $(\widehat{X}_i)_{1\leq i\leq p}$ be a continuous $p$-watermelon
without wall condition and $\widehat{\Sigma}_{k,p}$ be the $k$-th
elementary symmetric polynomial
of $(\widehat{X}_i)_{1\leq i\leq p}$, {\it i.e.}
$$\widehat{\Sigma}_{k,p}=\sum_{1\leq i_1<\dots<i_k\leq
p}\widehat{X}_{i_1}\dots \widehat{X}_{i_k}.$$
We have
\[
\esp{\widehat{\Sigma}_{2k,p}}=(-1)^k\frac{p!}{(p-2k)!2^kk!}\big(t(1-t)\big)^k
\qquad\text{and}\qquad
\esp{\widehat{\Sigma}_{2k+1,p}}=0.
\]
\end{prop}
\begin{proof}[Proof of Proposition \ref{moment avec mur}]
By Proposition \ref{densite}, we know that the density of a $p$-watermelon
with wall condition at time $t$ is $f(t;x)$, thus
\[
\esp{X_1(t)^k}=\int {x_1}^kf(t;x)\dd x
\qquad\text{and}\qquad
\esp{X_2(t)^k}=\int {x_2}^kf(t;x)\dd x.
\]
By carrying out the change of variables $u_i=x_i/\sqrt{t(1-t)}$, we obtain
\begin{multline*}
\esp{X_i(t)^k}=\frac{\big(t(1-t)\big)^{k/2}}{3\pi}
\int_{0\leq u_1\leq u_2}{u_i}^k({u_1}^2-{u_2}^2)^2{u_1}^2{u_2}^2\\
\times\exp\Big(-\frac{{u_1}^2+{u_2}^2}{2}\Big)\dd u_1\dd u_2.
\end{multline*}
Now, it suffices to compute this integral which we denote by
$I_k^i$. Expanding the integrated function and using integrations by parts,
we have
\[
I_k^1=-2\alpha_{k+5}+(16-k)\alpha_{k+3}+(k^2+2k+12)\int_{0\leq x\leq y}
x^{k+2}\,\textnormal{e}^{-\frac{x^2+y^2}{2}}\dd x\dd y
\]
where
\[
\alpha_n=\int_0^\infty x^n\,\textnormal{e}^{-x^2/2}\dd x
=\left\{
\begin{array}{ll}
\disp k!/2&\textnormal{ si }n=2k+1\\
\disp%2^{-2j}\binom{j}{2j}\sqrt{\pi}&\textnormal{ si }n=2k.
\frac{(2k)!}{2^k k!}\sqrt{\frac{\pi}{2}}&\textnormal{ si }n=2k.
\end{array}
\right.
\]
We compute this integral by new integrations by parts and we obtain
\begin{multline*}
I_{2k-1}^1=-2\alpha_{k+5}+(16-k)\alpha_{k+3}\\+(k^2+2k+12)
2^kk!\Big(\int_0^\infty\,\textnormal{e}^{-y^2/2}\dd y
-\sum_{j=0}^k2^{-j}\frac{1}{j!}\alpha_{2j}
\Big)
\end{multline*}
and
\begin{align*}
I_{2k}^1&=-2\alpha_{2k+5}+(16-2k)\alpha_{2k+3}\\
&\hspace{4.2cm}+4(k^2+k+3)
\int_{0\leq x\leq y}x^{2(k+1)}\,\textnormal{e}^{-(x^2+y^2)/2}\dd x\dd y\\
=&-2\alpha_{2k+5}+(16-2k)\alpha_{2k+3}\\
&\hspace{2.6cm}+4(k^2+k+3)
\frac{(2k+2)!}{2^{k+1}(k+1)!}
\Big(\frac{\pi}{4}
-\sum_{j=0}^k\alpha_{2j+1}\frac{(j+1)!2^{j+1}}{(2j+1)!}\Big).
\end{align*}
If we replace the $a_k$ by their values, we obtain the expected equality.

\smallskip

The computation of $I_{2k-1}^2$ is make by the same way.
\end{proof}

\begin{proof}[Proof of Proposition \ref{moment sans mur}]
As in the previous proof, we use the density $\widehat{f}$ of a
$2$-watermelon without wall condition.
\[
\esp{\widehat{X}_i(t)^k}
=\int_{x_1\leq x_2}{x_i}^k\widehat{f}(t;x)\dd x_1\dd x_2.
\]
Using the changes of variables $u_i=x_i/\sqrt{t(1-t)}$, we obtain
\[
\esp{\widehat{X}_i(t)^k}=\frac{1}{2\pi}\big(t(1-t)\big)^{k/2}
\int_{u_1\leq u_2}{u_i}^k(u_1-u_2)^2\exp(-\frac{{u_1}^2+{u_2}^2}{2})
\dd u_1\dd u_2.
\]
We denote by $I_k^i$ the above integral.

When $k$ is even, this integral is easily computable since, by symmetry,
\[
I_{2k}^1=I_{2k}^2
=\frac12\int_{\mathbb{R}^2}{x_1}^{2k}(x_1-x_2)^2
\exp(-\frac{{x_1}^2+{x_2}^2}{2})\dd x_1\dd x_2.
\]
By integrations by parts, we obtain
\[
I_{2k}^1=I_{2k}^2=
\pi\frac{(2k)!(k-1)}{2^{k+1}k!}.
\]

When $k$ in odd, we have
\[I_k^1=-I_k^2\]
and by integrations by parts
\begin{align*}
I_{2k}^1&=-2\int_{-\infty}^{+\infty}y^{2k+2}\,\textnormal{e}^{-y^2}\dd y
+(2k+3)\int_{x\leq y}x^{2k+1}\,\textnormal{e}^{-(x^2+y^2)/2}\dd x \dd y\\
&=-2\frac{(2k+2)!}{2^{2k+2}(k+1)!}\sqrt{\pi}
-\sum_{j=0}^k\frac{k!2^{j}}{j!}
\int_{-\infty}^{+\infty}x^{2j}\,\textnormal{e}^{-x^2}\dd x
\end{align*}
where the last integral is easy to compute.
\end{proof}

\begin{proof}[Proof of Proposition \ref{sym2}]
Let $\widehat{\Sigma}_{0,p}=1$, $\widehat{\Sigma}_{-1,p}=0$ and
$m_{k,p}(t)=\esp{\widehat{\Sigma}_{k,p}}$.
The It\^o's formula applied to $\widehat{\Sigma}_{k,p}$ yields
\begin{align*}
\dd\widehat{\Sigma}_{k,p}
&=\sum_{i=1}^p\Big(
\sum_{\substack{1\leq i_1<\dots<i_{k-1}\leq p \\ \forall
    l\in\{1,\dots,k-1\},\ i_l\neq i}}
\widehat{X}_{i_1}\dots \widehat{X}_{i_{k-1}}\Big)\dd \widehat{X}_i\\
&=\Big(\frac{-k\widehat{\Sigma}_{k,p}}{1-t}+S_{k,p}\Big)\dd
t+\widehat{\Sigma}_{k,p}\dd B,\\
\end{align*}
where
\[
S_{k,p}=\sum_{i=1}^p\Big(
\sum_{\substack{1\leq i_1<\dots<i_{k-1}\leq p \\ \forall
    l\in\{1,\dots,k-1\},\ i_l\neq i}}
\widehat{X}_{i_1}\dots \widehat{X}_{i_{k-1}}\sum_{\substack{1\leq j\leq p\\j\neq i}}
\frac{1}{\widehat{X}_i-\widehat{X}_j}\Big)
\]
and $\widehat{\Sigma}_{k,p}^2$ is a polynomial of $(\widehat{X}_i)$.
We remark that, in $S_{k,p}$, the terms, where all $i_l$ are not equal to
$j$, cancel out.
\begin{align*}
S_{k,p}
&=\sum_{1\leq i\neq j\leq p}\Big(
\sum_{\substack{1\leq i_1<\dots<i_{k-1}\leq p \\ \forall
    l\in\{1,\dots,k-2\},\ i_l\neq i,j}}
\frac{\widehat{X}_j\widehat{X}_{i_1}\dots \widehat{X}_{i_{k-2}}}{\widehat{X}_i-\widehat{X}_j}\Big)\\
&=\frac12\sum_{1\leq i\neq j\leq p}\Big(
\sum_{\substack{1\leq i_1<\dots<i_{k-1}\leq p \\ \forall
    l\in\{1,\dots,k-2\},\ i_l\neq i,j}}\widehat{X}_{i_1}\dots \widehat{X}_{i_{k-2}}
\big(\frac{\widehat{X}_j}{\widehat{X}_i-\widehat{X}_j}+\frac{\widehat{X}_i}{\widehat{X}_j-\widehat{X}_i}\big)\Big)\\
&=-\frac12\sum_{1\leq i\neq j\leq p}\Big(
\sum_{\substack{1\leq i_1<\dots<i_{k-1}\leq p \\ \forall
    l\in\{1,\dots,k-2\},\ i_l\neq i,j}}\widehat{X}_{i_1}\dots \widehat{X}_{i_{k-2}}
\Big)\\
&=-\frac{\big(p-(k-2)\big)\big(p-(k-1)\big)}{2}\widehat{\Sigma}_{k-2,p}.
\end{align*}
Hence, $\widehat{\Sigma}_{k,p}$ satisfies the SDE
\[\dd\widehat{\Sigma}_{k,p}=\Big(-\frac{k\widehat{\Sigma}_{k,p}}{1-t}
-\frac{\big(p-(k-2)\big)\big(p-(k-1)\big)}{2}\widehat{\Sigma}_{k-2,p
}\Big)\dd t
+\widehat{\Sigma}_{k,p}\dd B.
\]
When we take the expectation in the both side of this SDE, we obtain
\[m'_{k,p}(t)=\frac{km_{k,p}(t)}{1-t}
-\frac{\big(p-(k-2)\big)\big(p-(k-1)\big)}{2}m_{k-2,p}(t).\]
We know moreover that
\[m_{k,p}(t)=\int\widehat{\Sigma}_{k,p}(x_1,\dots,x_p)
\widehat{f}(t;x_1,\dots,x_p)\dd x_1\cdots \dd x_p.\]
The change of variable $\sqrt{t(1-t)}u_i=x_i$ yields
$$m_{k,p}(t)=c_{k,p}\big(t(1-t)\big)^{k/2},$$
where $c_{k,p}$ is a constant. It is clear that $c_{2k+1,p}$ is not
equal to zero. When $m_{2k,p}$ is replaced by $c_{2k,p}\big(t(1-t)\big)^{k}$ in
the SDE, we show that $c_{2k,p}$ satisfies the following recurrence
relation
$$c_{2k,p}=-\frac{\big(p-(2k-2)\big)\big(p-(2k-1)\big)}{2k}c_{2k-2,p}.$$
Since $c_{0,p}=1$, it comes
$$c_{2k,p}=(-1)^k\frac{p!}{(p-2k)!2^kk!},$$
this complete the proof of Proposition \ref{sym2}.
\end{proof}

\begin{proof}[Proof of Proposition \ref{sym1}]
Let $\Sigma_{0,p}=1$. Using the It\^o's formula, we show that $\Sigma_{k,p}$
satisfies the SDE
$$\dd \Sigma_{k,p}=
\Big(\frac{-2k\Sigma_{k,p}}{1-t}+(p-(k-1))(2(p-k)+3)\Sigma_{k-1,p}\Big)
\dd t+\Sigma_{k,p}\dd B,$$
where ${\Sigma_{k,p}}^2$ is a polynomial of $(X_i)$. Taking the expectation
in the above SDE, we obtain a partial differential equation satisfy by
$m_{k,p}$:
$$m_{k,p}'(t)=\frac{-2km_{k,p}(t)}{1-t}+(p-(k-1))(2(p-k)+3)m_{k-1,p}(t).$$
Since $f(t;.)$ is the density if the random variable $X(t)$, we show that
$m_{k,p}(t)=c_{k,p}\big(t(1-t)\big)^k$.
The previous partial differential equation give us the following recurrence
relation satisfied by $c_{k,p}$:
$$c_{k,p}=\frac{(p-(k-1))(2(p-k)+3)}{k}c_{k-1,p}.$$
Hence, we have finally
$$c_{k,p}=\frac{(2p+1)!}{(2p+1-2k)!2^kk!}.$$
\end{proof}

\bigskip

%%%%%%%%%%%%%%%%%%%%%%%%%%%%%%%%%%%%%%%%%%%%%%%%%%%%%%%%%%%%%%%%%%%%%%%%%%%%
%%%%%%%%%%%%%%%%%%%%%%%%%%%%%%%%%%%%%%%%%%%%%%%%%%%%%%%%%%%%%%%%%%%%%%%%%%%%

%\nocite{bonichon,bonichon2}

\bibliographystyle{plain}
%\bibliography{biblio}
%\end{document}

\end{document}